\documentclass[12pt,a4paper]{paper}

\usepackage{amssymb,amsmath, amsthm, amsfonts}
\usepackage{enumerate}
\usepackage{geometry}
\usepackage{tikz}
\usepackage{footmisc}

\newcommand{\termvalue}[1] {\lbrack\!\lbrack #1 \rbrack\!\rbrack}

\newcommand{\axD}{\textsf{(D)}}
\newcommand{\axK}{\textsf{(K)}}
\newcommand{\axKu}{\textsf{(K1)}}
\newcommand{\axKd}{\textsf{(K2)}}
\newcommand{\axMu}{\textsf{(M1)}}
\newcommand{\axMd}{\textsf{(M2)}}
\newcommand{\axMt}{\textsf{(M3)}}
\newcommand{\axMc}{\textsf{(M4)}}
\newcommand{\axDNu}{\textsf{(DN1)}}
\newcommand{\axDNd}{\textsf{(DN2)}}
\newcommand{\axT}{\textsf{(T)}}
\newcommand{\axBF}{\textsf{(BF)}}
\newcommand{\axCBF}{\textsf{(CBF)}}
\newcommand{\axNBF}{\textsf{(NBF)}}
\newcommand{\axPBF}{\textsf{(PBF)}}
\newcommand{\axNeq}{\textsf{(N$=$)}}
\newcommand{\axPeq}{\textsf{(P$=$)}}
\newcommand{\axu}{\textsf{(Ax1)}}
\newcommand{\axd}{\textsf{(Ax2)}}
\newcommand{\axt}{\textsf{(Ax3)}}
\newcommand{\axq}{\textsf{(Ax4)}}
\newcommand{\axc}{\textsf{(Ax5)}}
\newcommand{\axs}{\textsf{(Ax6)}}
\newcommand{\axst}{\textsf{(Ax7)}}
\newcommand{\axo}{\textsf{(Ax8)}}
\newcommand{\axCequ}{\textsf{(C$=$1)}}
\newcommand{\axCeqd}{\textsf{(C$=$2)}}
\newcommand{\axfour}{\textsf{(4)}}
\newcommand{\axfive}{\textsf{(5)}}

\newcommand{\axBFn}{\textsf{(BF$_n$)}}
\newcommand{\axCBFn}{\textsf{(CBF$_n$)}}
\newcommand{\axBFd}{\textsf{(BF$_2$)}}
\newcommand{\axCBFd}{\textsf{(CBF$_2$)}}
\newcommand{\axBFe}{\textsf{(BF$_\strictif$)}}
\newcommand{\axCBFe}{\textsf{(CBF$_\strictif$)}}

\DeclareSymbolFont{symbolsC}{U}{txsyc}{m}{n}
\DeclareMathSymbol{\strictif}{\mathrel}{symbolsC}{74}

\newtheorem{definition}{\vspace{1mm}Definition}[section]
\newtheorem{remark}[definition]{\vspace{1mm}Remark}
\newtheorem{lemma}[definition]{\vspace{1mm}Lemma}
\newtheorem{theorem}[definition]{\vspace{1mm}Theorem}

\newtheorem{proposition}[definition]{\vspace{1mm}Proposition}

\title{Modal Logic With Non-deterministic Semantics: Part~II - Quantified Case}
\author{Marcelo E.~Coniglio$^{1}$, Luis Fari\~nas del Cerro$^{2}$ and Newton M. Peron$^{3}$\\ [2mm] %
{\small $^1$Institute of Philosophy and the Humanities (IFCH) and}\\[-2mm]
{\small Centre for Logic, Epistemology and The History of Science (CLE),}\\[-2mm]
{\small University of Campinas (UNICAMP), Campinas, SP, Brazil.}\\[-2mm]
{\small E-mail:~{coniglio@unicamp.br}}\\[1mm]
{\small $^2$IRIT, Université de Toulouse,  CNRS, France. }\\[-2mm]
{\small E-mail:~{farinas@irit.fr}}\\[1mm]
{\small $^3$Graduate Program in Philosophy,}\\[-2mm]
{\small Federal University of Southern Frontier (UFFS), Chapec\'o, SC, Brazil.}\\[-2mm]
{\small E-mail:~{newton.peron@uffs.edu.br}}}

\begin{document}

\maketitle

\begin{abstract}
In the first part of this paper we  analyzed  finite non-deterministic matrix semantics for propositional non-normal modal logics as an alternative to the standard  Kripke's possible world  semantics. This kind of modal systems characterized by finite non-deterministic matrices was originally proposed by Ju.~Ivlev in the 70's. The aim of this second paper is to introduce a formal non-deterministic semantical framework for the quantified versions of some Ivlev-like non-normal modal logics. It will be shown that several well-known controversial issues of quantified modal logics, relative to the identity predicate, Barcan's formulas, and de dicto and de re modalities,  can be tackled from a new angle  within the present framework.
%can be avoided or better controlled within the present framework.      
\end{abstract}

\


\section*{Introduction}

In previous papers (see~\cite{con:cer:per:15,con:cer:per:17,con:cer:per:19}),\footnote{Closely related results were independently obtained by Omori and Skurt in~\cite{omo:sku:16}.} we  analyzed  finite non-deterministic matrix semantics for propositional modal logics as an alternative to the standard  Kripke's possible world  semantics, in order to better understand the modal concepts of `necessary' and `possible' at the propositional level. This kind of semantics was independently proposed by Y. Ivlev (see~\cite{ivl:73,ivl:85,ivl:88,ivl:13}) and  J. Kearns (see~\cite{kear:81}). The aim of this paper, which is the second part of~\cite{con:cer:per:19}, is to introduce a formal non-deterministic semantical framework for some non-normal first-order modal logics, based on Ivlev's  approach to propositional non-normal modal logics.
In such works,  Ivlev introduced  several modal systems which do not have the necessitation rule, for instance, weaker versions of {\bf T} and {\bf S5}. The semantics proposed by Ivlev is given by finite-valued  Nmatrices together with the  notion of (single-valued) valuations considered by A. Avron and I. Lev in~\cite{avr:lev:01} (which also introduced the terminology {\em non-deterministic matrices/Nmatrices} and  {\em legal valuations}). Thus, Ivlev's systems constitute one of the earliest antecedents of Avron and Lev's non-deterministic semantics.\footnote{Other antecedents of Nmatrices were proposed by N. Rescher~\cite{res:62}, Quine~\cite{qui:74}, Kearns~\cite{kear:81},  and J. Crawford  and D. Etherington~\cite{craw:ethe:98}.} The quantified version of the Ivlev-like modal logics presented here is along the same lines as the non-deterministic semantics proposed in~\cite{CFG19}  for quantified paraconsistent logics. As we shall see along this paper, several well-known criticisms to quantified modal logics, relative to the identity predicate, Barcan's formulas, and de dicto and de re modalities,  can be tackled from a new angle  within the present framework.   
	
In Part~I of this paper (\cite{con:cer:per:19}), Ivlev's propositional systems were expanded  and it was shown, among other results, that several Ivlev-like modal systems, which are characterized by  finite-valued non-deterministic semantics, can be captured in terms of the modal concepts of \emph{necessarily true}, \emph{possibly true} and \emph{actually true}.\footnote{This is closely related to the swap-structures semantics proposed in~\cite{con:gol:19} for some of these systems.} In particular, the four-valued systems  can be captured  by only two concepts:  \emph{actually true} and \emph{contingently true}. 

Now, let us consider the quantified case. In first-order classical logic, a  (usually called {\em Tarskian}) structure is a pair $\mathfrak{A} = \langle U, \cdot^\mathfrak{A} \rangle$ such that $U$ is a nonempty set called the \emph{domain} or {\em universe} of the structure, in which the individuals of the structure exist. The function $\cdot^\mathfrak{A}$ assigns a concrete interpretation for the symbols of the signature. In particular,  $n$-ary  predicate symbols are interpreted as set of $n$-tuples over $U$, which is the extension of the predicate in $\mathfrak{A}$.\footnote{See, for instance, \cite[p.  49]{men:10} and~\cite[p. 80--81]{ender:01}.}

In Kripke semantics for first-order modal logic, there are at least two universes for that structure. First, there is a set of \emph{worlds} and a relation between them. Besides this, the function $\mathfrak{A}$ associates to each predicate its extension in a given world. In a {\em constant domain} approach, there is a fixed nonempty set $U$ which states the individuals that exist in every world. In a {\em varying domain} approach, the set $U$ is replaced by a function $U(w)$ that stablishes the domain of the individuals that exist in the \emph{world} $w$.\footnote{For constant domains, see \cite[p. 243]{hug:cre:96} and \cite[p. 95-98]{fit:med:98}. For varying domains, see \cite[p. 278--279]{hug:cre:96} and \cite[p. 101--104]{fit:med:98}.}

%In Part~I of this paper (\cite{con:cer:per:19}), Ivlev systems were expanded  and it was shown, among other results, that several Ivlev-like modal systems, which are characterized by  finite-valued non-deterministic semantics, can be captured in terms of the modal concepts of \emph{necessarily true}, \emph{possibly true} and \emph{actually true}.\footnote{This is closely related to the swap-structures semantics proposed in~\cite{con:gol:19} for some of these systems.} In paricular, the four-valued systems  can be captured  by only two concepts:  \emph{actually true} and \emph{contingently true}.

Let {\bf Tm} be the Ivlev's four-valued version of {\bf T}, and let us denote by {\bf Tm$^*$} the proposed quantified version of it.
Since the four-valued non-deterministic semantics for {\bf Tm} can be described in terms of the concepts of \emph{actually true} and \emph{contingently true}, this suggest that  the predicate symbols in a first-order structure $\mathfrak{A}$ for {\bf Tm$^*$} could be interpreted in terms of two mappings $\mathfrak{a}^\mathfrak{A}$ and $\mathfrak{c}^\mathfrak{A}$, describing the actual and the contingent extension of $P$ in $\mathfrak{A}$. In more general systems involving six or eight truth-values, which are explained in terms of the modal concepts of \emph{necessarily true}, \emph{possibly true} and \emph{actually true}, each predicate could be interpreted in terms of three mappings describing the respective extensions.

The organization of this paper is as follows: Section~\ref{intuition} presents the semantical intuitions behind the non-deterministic framework for the modal systems proposed here. In Section~\ref{Tm*}  the basic four-valued system {\bf Tm$^*$} is introduced, as a possible first-order version of Ivlev's four-valued system {\bf Tm}. As we shall see, at least two possible (non-deterministic) interpretations for the quantifiers could be considered. Section~\ref{filosofia} analyzes some well-known philosophical questions related to the standard approach to quantified modal logic, showing that some of the criticisms to quantified modal logics can be can be avoided or better controlled  in  {\bf Tm$^*$}. In particular, it will be shown that the choice of one or another interpretation for the quantifiers in   {\bf Tm$^*$} is related to the distinction between {\em de dicto} and {\em de re} modalities. In Section~\ref{other-sys}, first-order versions of other four-valued, six-valued and eight-valued Ivlev-like modal systems are briefly discussed. Finally, Section~\ref{final} discusses the results presented along the paper as well as some future lines of research. 
 
%The key-idea of this Part~II is to extend standard first-order structures to four-valued non-deterministic structures endowed with two interpretation functions for the predicate symbols, each function capturing each of those two modal concepts. The generalization to first-order logics based on some six-valued and eight-valued systems studied in Part~I will be briefly discussed in Section~\ref{other-sys}.

\section{The semantical intuitions behind (first-order) non-deterministic modal semantics} \label{intuition} 
%\section{Semantic intuition} \label{intuition}

According to the brief discussion presented in the Introduction, the four-valued non-deterministic semantics for {\bf Tm} can be described in terms of the concepts of \emph{actually true} and \emph{contingently true}, hence the predicate symbols in a structure $\mathfrak{A}$ for {\bf Tm$^*$} can be interpreted in terms of two mappings $\mathfrak{a}^\mathfrak{A}$ and $\mathfrak{c}^\mathfrak{A}$. To this end, in this paper we will consider semantical structures for a non-deterministic four-valued first-order modal logic as being pairs $\mathfrak{A}=\langle U,  \cdot^\mathfrak{A}\rangle$ such that the interpretation under $\cdot^\mathfrak{A}$  of function symbols and individual constants is  defined as usual in Tarskian first-order structures. On the other hand, each $n$-ary predicate symbol $P$ will be  interpreted as a pair $P^\mathfrak{A}(P)=(\mathfrak{a}^\mathfrak{A}(P), \mathfrak{c}^\mathfrak{A}(P))$ of subsets of $U^n$.\footnote{This idea was inspired on the first-order extension of Kleene's three-valued logic K3 proposed by  Kripke in~\cite{kri:75}.} The function $\mathfrak{a}^\mathfrak{A}$ says whenever a tuple of  individuals of the domain  {\em actually} satisfies  or not predicate $P$. In other words, the function $\mathfrak{a}^\mathfrak{A}$ assigns to the predicate $P$ its {\em actual} extension in $\mathfrak{A}$. In turn, the function $\mathfrak{c}^\mathfrak{A}$ assigns to $P$ its  \emph{contingent} extension in $\mathfrak{A}$. It says whenever a tuple of  individuals of the domain {\em  contingently} satisfies  or not predicate $P$.
These two modal concepts produce four modal values that can be interpreted as follows:

\begin{itemize}
  \item[] \  $T^+$: necessarily true;
  \item[] \  $C^+$: contingently true;
  \item[] \  $C^-$: contingently false;
  \item[] \  $F^-$: necessarily false / impossible.
\end{itemize}

Thus, the interpretation of a given $n$-ary predicate $P$ in a four-valued structure $\mathfrak{A}$ as above gives origin to the following configuration over  $U^n$:

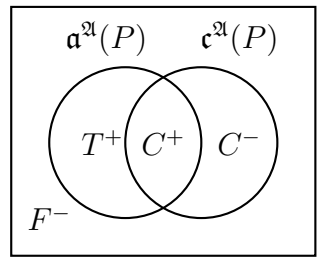
\begin{center}
\begin{tikzpicture}[thick] 
\draw (2.5,-1.5) rectangle (-1.5,1.8) node[below right] {}; 
\draw (0,0) circle (1) node[above,shift={(0,1)}] {\hspace*{-5mm}$\mathfrak{a}^\mathfrak{A}(P)$};
\draw (1,0) circle (1) node[above,shift={(0,1)}] {\hspace*{10mm}$\mathfrak{c}^\mathfrak{A}(P)$};
	\node at (-0.3,0) {$T^+$}; 
	\node at (0.5,0) {$C^+$}; 
	\node at (1.5,0) {$C^-$};	
	\node at (-1,-1) {$F^-$}; 
\end{tikzpicture}
\end{center}
 
The truth-value attached to each of the four areas above is the value assigned to the atomic formula $P(\tau_1,\ldots,\tau_n)$ when the $n$-tuple  $(\tau_1^\mathfrak{A},\ldots,\tau_n^\mathfrak{A})$ of $U^n$ associated to $(\tau_1,\ldots,\tau_n)$ in  $\mathfrak{A}$   belongs to that area. For instance, $P(\tau_1,\ldots,\tau_n)$ gets the value $T^+$ in $\mathfrak{A}$ iff $(\tau_1^\mathfrak{A},\ldots,\tau_n^\mathfrak{A}) \in \mathfrak{a}^\mathfrak{A}(P) \setminus \mathfrak{c}^\mathfrak{A}(P)$. In  Remark~\ref{equiv-sem} it will be shown that this approach is equivalent to consider, for every $n$-ary predicate symbol $P$, a function $P_\mathfrak{A}:U^n \to V_4$, where $V_4$ is the set of four truth-values corresponding to the four areas in the figure above.\footnote{Observe that the interpretation of predicate symbols as mappings from the domain of a first-order structure into the domain of an algebra corresponds to what is  usually done in the realm of algebraic semantics for first-order logics.}

In order to better understand  the first-order expansion of this semantics in a intuitive way, let us consider the following sentences as example:

\begin{center}
(1) \ \ \  Socrates is mortal
\end{center}

The sentence (1) is  necessarily true when the individual Socrates is in the actual extension of the predicate 
``mortal'' but not in the contingent extension of it. But (1) is only contingently true when Socrates is both in the actual and in the contingent extension of the predicate ``mortal''. If  Socrates is not in the actual but is in the contingent extension of the predicate ``mortal'', we say that (1) is contingently false. Finally, if  Socrates is neither in the actual nor in the contingent extension of the predicate ``mortal'', we conclude that (1) is impossible.

Consider now the universal sentence:

\begin{center}
 (2) \ \ \  All are  mortals
\end{center}

The sentence (2) will be necessarily true when it is necessarily true for each individual of the domain. In this case, the actual extension of the predicate ``mortal'' must coincide with the domain. Besides, the contingent extension of ``mortal'' must be empty. The sentence (2)  is only contingently true if: (i) there is at least one individual of the domain that has contingently the propriety of being mortal; and (ii) each individual of the domain is necessarily mortal or only contingently mortal. In order to guarantee both conditions, we require that the actual extension of ``mortal'' coincides with the domain, but its contingent extension must not be empty. 

We say that (2) is contingently false if at least one individual of the domain has not the propriety of being mortal. Besides, any individual of the domain \emph{could} be mortal, that is, they have contingently the propriety of being mortal. The idea here is that the actual extension of ``mortal'' does not coincide with the domain but the union of the actual and the contingent extension of ``mortal'' does. 
Finally, (2)  is impossible when it is not possible to at least one individual to have the propriety of being mortal, that is, at least one thing is not actually neither contingently mortal. We can state that condition saying that there is at least one individual of the domain that is not in the union of the actual and the contingent extension of ``mortal''.

With respect to equality,  consider the sentence:

\begin{center}
(3) \ \ \   Phosphorus is Hesperus
\end{center}

In Kripke semantics (with constant or varying domains), proper names refer to the same individual through all the possible worlds. Because of this, they are called \emph{rigid designators}. If two proper names refer to the same individual of the domain in some possible world, they will refer to the same one through all the possible worlds, that is, equalities are always necessary.\footnote{See \cite[p. 48]{kri:81}. The sentence (3) is an  example given by Kripke in~\cite[p. 28--29]{kri:81}. } Conversely --- with certain restrictions on the accessibility relation between the worlds in our Kripkean models --- if two proper names refer to different individuals in some possible world, they will refer to different ones through all the possible worlds. That means, inequalities are necessary too.\footnote{In \cite[p. 311]{hug:cre:96}, the authors showed also that in Kripkean models in which the relation between possible world is, at least, reflexive and symmetrical, the following holds: $$(x\not\approx y) \to  \Box (x \not\approx y).$$ This proves that in some Kripkean systems not only equalities of any kind (between constants or functions) are necessary but also that inequalities of any kind are necessary\label{converse_rigidity}.}

Those considerations force us to admit that, in Kripke semantics, (3) is either necessarily true or necessary false, that is, impossible. In this section, we will present an analogous approach considering first-order version of Ivlev semantic, that is, (3) will  receive only two values: the value ``necessarily true'' when ``Phosphorus'' and ``Hesperus'' refer to the same individual in the domain; the value ``impossible'' if they don't. However, this can be seen as an arbitrary option. In Subsection~\ref{contingent}, we will discuss how Ivlev's semantics can also deal with contingent identities in a  natural way.

\section{The first-order non-normal modal logic {\bf Tm}$^*$} \label{Tm*}

In this section the intuitive ideas described in Section~\ref{intuition} will be formalized. Recall that {\bf Tm} is one of the (non-normal) propositional modal systems introduced by Ivlev in~\cite{ivl:88}, but also studied  independently by Kearns in~\cite{kear:81}, and more recently by Omori and Skurt in~\cite{omo:sku:16} and by us in~\cite{con:cer:per:15, con:cer:per:17, con:cer:per:19}. An outstading feature of {\bf Tm} is that it is semantically characterized by a four-valued non-deterministic matrix (or Nmatrix), by using the terminology and formalism introduced by Avron and Lev in~\cite{avr:lev:01} (see also~\cite{avr:lev:05}). Thus, Ivlev's modal logics constitute an early antecendent  of Nmatrices (called {\em quasi matrices} by Ivlev in~\cite{ivl:88}). 

\subsection{The logic {\bf Tm} and  two ways of defining quantifiers over it} \label{Nmatrix}

The propositional modal logic {\bf Tm} is defined over the propositional signature  $\Sigma$, which consists of the connectives $\neg$ (negation), $\Box$ (necessary) and $\to$ (implication). Semantically,  {\bf Tm} is
characterized by a four-valued  non-deterministic matrix  $\mathcal{M}_{\bf Tm} = \langle \mathcal{A}_{\bf Tm},D\rangle$ such that $\mathcal{A}_{\bf Tm} = \langle V_4,\cdot^{\bf Tm}\rangle$ is a multialgebra over $\Sigma$ (recall~\cite[Definition~2.1]{con:cer:per:19}) with domain $V_4=\{T^+,C^+,C^-,F^-\}$ and $D=\{T^+,C^+\}$ is the set of {\em designated} truth-values.   From now on the multioperation $\#^{\bf Tm}$ associated to each connective $\#$ will be simply denoted by $\tilde{\#}$, when there is no risk  of confusion. The connectives of  $\Sigma$ are interpreted in $\mathcal{A}_{\bf Tm}$ as follows:

	\begin{displaymath}
		\begin{array}{c c c c}
			\begin{array}{|c|c|}
				\hline x & \tilde{\neg} x\\
				\hline T^+ & \{F^-\}  \\
				\hline C^+ & \{C^-\}  \\
				\hline C^- & \{C^+\} \\
				\hline F^- & \{T^+\} \\
				\hline
			\end{array}
			&
			\begin{array}{|c|c|}
				\hline x & \tilde{\Box} x\\
				\hline T^+ & \{T^+,C^+\}  \\
				\hline C^+ & \{C^-,F^-\}  \\
				\hline C^- & \{C^-,F^-\} \\
				\hline F^- & \{C^-,F^-\} \\
				\hline
			\end{array}
			&
			
			\begin{array}{|c|c|c|c|c|}
				\hline \tilde{\to}  	& T^+ & C^+	& C^-	& F^- \\
				\hline 											T^+	 	& \{T^+\} & \{C^+\} 	& \{C^-\} 	& \{F^-\}\\
				\hline 											C^+	 	& \{T^+\} & \{T^+,C^+\} & \{C^-\} 	& \{C^-\}\\
				\hline 											C^-	 	& \{T^+\} & \{T^+,C^+\} & \{T^+,C^+\} 	& \{C^+\}\\
				\hline 											F^-		& \{T^+\} & \{T^+\} 	& \{T^+\} 	& \{T^+\}\\
				\hline
			\end{array}
		\end{array}
	\end{displaymath}

Disjunction and conjunction can be defined in {\bf Tm} as follows: $\alpha \vee \beta := \neg \alpha \to \beta$ and $\alpha \wedge \beta := \neg(\alpha \to \neg\beta)$, while possibility is given as usual by $\Diamond \alpha:=\neg\Box\neg\alpha$. The corresponding multioperators are defined in $\mathcal{A}_{\bf Tm}$ as follows:

\begin{displaymath}
\begin{array}{|c|c|}
				\hline x & \tilde{\Diamond} x\\
				\hline T^+ & \{T^+,C^+\}  \\
				\hline C^+ & \{T^+,C^+\}  \\
				\hline C^- & \{T^+,C^+\} \\
				\hline F^- & \{C^-,F^-\} \\
				\hline
			\end{array}
\hspace*{1cm}
\begin{array}{|c|c|c|c|c|}
      \hline \tilde{\vee} 	& T^+ 	& C^+ & C^- & F^-\\
\hline \hline
T^+ & \{T^+\} & \{T^+\}  & \{T^+\} &  \{T^+\}\\
       \hline
C^+ & \{T^+\} & \{T^+,C^+\}  & \{T^+,C^+\} 		& \{C^+\}\\
       \hline
C^- & \{T^+\} & \{T^+,C^+\} & \{C^-\} & \{C^-\}\\
       \hline
F^- & \{T^+\} & \{C^+\} & \{C^-\} & \{F^-\}\\    
       \hline
\end{array}
\end{displaymath}

\

\begin{displaymath}
\begin{array}{|c|c|c|c|c|}
      \hline \tilde{\wedge} 	& T^+ 	& C^+ & C^- & F^-\\
\hline \hline
T^+ & \{T^+\} & \{C^+\}  & \{C^-\} &  \{F^-\}\\
       \hline
C^+ & \{C^+\} & \{C^+\}  & \{F^-,C^-\} 	& \{F^-\}\\
       \hline
C^- & \{C^-\} & \{F^-,C^-\} & \{F^-,C^-\} & \{F^-\}\\
       \hline
F^- & \{F^-\} & \{F^-\} & \{F^-\} & \{F^-\}\\    
       \hline
\end{array}
\end{displaymath}

\

It is easy to see that $x \,\tilde{\vee}\, y = \tilde{\neg}(\tilde{\neg} x \,\tilde{\wedge}\, \tilde{\neg} y)$ and  $x \,\tilde{\wedge}\, y = \tilde{\neg}(\tilde{\neg} x \,\tilde{\vee}\, \tilde{\neg} y)$ for every $x,y \in V_4$. That is, the De Morgan rules are valid in $\mathcal{A}_{\bf Tm}$.

In Subsection~\ref{axioTm*} a  first-order extension of {\bf Tm} called  {\bf Tm}$^*$ will be proposed.\footnote{The first-order extension of other Ivlev-like modal systems as the ones studied in the first part of this paper~\cite{con:cer:per:19} can be done in an analogous way, by means of straightforward adaptations.} The semantics of {\bf Tm}$^*$ will be given by means of first-order structures $\mathfrak{A}$ evaluated over the Nmatrix $\mathcal{M}_{\bf Tm}$. This approach to first-order logics with a non-deterministic semantics based on Nmatrices was already considered in the literature for several paraconsistent logics known as {\em logics of formal inconsistency},  see~\cite{avr:zam:07} and~\cite{CFG19}. The latter considers a family of Nmatrices defined over certain non-deterministic algebras called {\em swap structures}. These structures, which generalize the finite-valued Nmatrices, were developed in
~\cite{con:gol:19} for several Ivlev-like modal logics. As we shall see, the semantics for {\bf Tm}$^*$ adapted from~\cite{avr:zam:07} and~\cite{CFG19} will be equivalent to the  semantical approach informally described in the previous section. 

The first step is extending the Nmatrix $\mathcal{M}_{\bf Tm}$ with a multioperator $\tilde{Q}_4: (\mathcal{P}(V_4)-\{\emptyset\}) \to (\mathcal{P}(V_4)-\{\emptyset\})$ for every quantifier $Q\in \{\forall, \exists\}$. The idea is that a given valuation over the extended Nmatrix will choose  a value, for a given formula of the form $Q x \varphi$, within the set $\tilde{Q}_4(X)$, where $X$ is the set of instances of $\varphi(x)$ over the given first-order structure $\mathfrak{A}$.\footnote{For `instances of $\varphi(x)$ over $\mathfrak{A}$' we mean  the set of denotations of $\varphi$ over the structure $\mathfrak{A}$, when $x$ takes all the possible values in the domain of $\mathfrak{A}$, assuming that any variable occurring free in $\varphi$ other than $x$ has ben interpreted by a given assignment over $\mathfrak{A}$. The technical details will be given in Subsection~\ref{structures}.} Accordingly to the previous approaches to quantified Nmatrices, which are inspired in algebraic quantified logics, it is natural that $\tilde{\forall}_4(X)$ and $\tilde{\exists}_4(X)$ be defined as the conjunction  and the disjunction of the members of $X$ according to the respective multioperators of  $\mathcal{A}_{\bf Tm}$. This produces the following:

	\begin{displaymath}
			\begin{array}{|c|c|}
				\hline  X & \tilde{\forall}_4(X)\\
				\hline \{T^+\} & \{T^+\}  \\
				\hline \{C^+\} & \{C^+\}  \\
				\hline \{T^+, C^+\} & \{C^+\} \\
				\hline \{C^-, C^+\} & \{F^-, C^-\} \\
				\hline \{C^-, C^+, T^+\} & \{F^-, C^-\} \\
				\hline \{C^-\} & \{C^-\} \\
				\hline \{C^-, T^+\} & \{C^-\} \\
				\hline F^- \in X & \{F^-\} \\
				\hline
			\end{array}
			\hspace{1cm}
			\begin{array}{|c|c|}
				\hline  X & \tilde{\exists}_4(X)\\
				\hline T^+ \in X & \{T^+\}  \\
				\hline \{C^+\} & \{C^+\}  \\
				\hline \{C^+, F^-\} & \{C^+\} \\
				\hline \{C^+, C^-\} & \{T^+, C^+\} \\
				\hline \{C^+, C^-, F^-\} & \{T^+, C^+\} \\
				\hline \{C^-\} & \{C^-\} \\
				\hline \{C^-, F^-\} & \{C^-\} \\
				\hline \{F^-\} & \{F^-\} \\
				\hline
			\end{array}
	\end{displaymath}

It is easy to see that $\tilde{\forall}_4(X) = \tilde{\neg} \tilde{\exists}_4(\tilde{\neg} X)$ and $\tilde{\exists}_4(X) = \tilde{\neg} \tilde{\forall}_4(\tilde{\neg} X)$, where, for every $\emptyset\neq X \subseteq V_4$, $\tilde{\neg} X=\{\tilde{\neg} x \ : \ x \in X\}$. However, since we are interested in  the satisfaction of the Barcan formulas (see Subsection~\ref{sectbarcan}), and taking into consideration the semantic intuitions of the sentence (2) explored in section \ref{intuition}, stricter (deterministic) forms of quantification will be considered in {\bf Tm}$^*$, namely $\tilde{\forall}_4^d$ and $\tilde{\exists}_4^d$, given by the tables above.

	\begin{displaymath}
			\begin{array}{|c|c|}
				\hline  X & \tilde{\forall}_4^d(X)\\
				\hline \{T^+\} & \{T^+\}  \\
				\hline \{C^+\} & \{C^+\}  \\
				\hline \{T^+, C^+\} & \{C^+\} \\
				\hline \{C^-, C^+\} & \{C^-\} \\
				\hline \{C^-, C^+, T^+\} & \{C^-\} \\
				\hline \{C^-\} & \{C^-\} \\
				\hline \{C^-, T^+\} & \{C^-\} \\
				\hline F^- \in X & \{F^-\} \\
				\hline
			\end{array}
			\hspace{1cm}
			\begin{array}{|c|c|}
				\hline  X & \tilde{\exists}_4^d(X)\\
				\hline T^+ \in X & \{T^+\}  \\
				\hline \{C^+\} & \{C^+\}  \\
				\hline \{C^+, F^-\} & \{C^+\} \\
				\hline \{C^+, C^-\} & \{C^+\} \\
				\hline \{C^+, C^-, F^-\} & \{C^+\} \\
				\hline \{C^-\} & \{C^-\} \\
				\hline \{C^-, F^-\} & \{C^-\} \\
				\hline \{F^-\} & \{F^-\} \\
				\hline
			\end{array}
	\end{displaymath}

It is easy to see that the deterministic quantifiers correspond, respectively, to the deterministic conjunction  and disjunction of the members of $X$ according to the order given by the chain $F^- \leq C^- \leq C^+ \leq T^+$.
Clearly $\tilde{\forall}_4^d(X) = \neg \tilde{\exists}_4^d(\neg X)$ and $\tilde{\exists}_4^d(X) = \neg \tilde{\forall}_4^d(\neg X)$ for every $\emptyset\neq X \subseteq V_4$.

As we shall see  in Subsection~\ref{de_re}, in order to analyze the distinction between  the {\em de re} and {\em de dicto} modalities, both kinds of quantifiers $\tilde{Q}_4$ and $\tilde{Q}_4^d$ will be relevant.

\subsection{The logic {\bf Tm}$^*$ and its axiomatics} \label{axioTm*}

To our purposes, we will consider first-order modal languages based on the propositional signature $\Sigma$ described at the beginning of Subsection~\ref{Nmatrix}, expanded with the universal quantifier $\forall$,\footnote{As discussed in the previous subsection, conjunction $\wedge$, disjunction $\vee$ and existential quantifier $\exists$ can be defined in terms of the other symbols, hence it will be ommited from the list of primitive symbols.} and defined over first-order signatures, which are defined as usual. Thus, a {\em first-order signature} is a collection $\Theta$ formed by the following symbols:  (i) a non-empty set of predicate symbols $\mathcal{P}$, with the corresponding arity $\varrho(P) \geq 1$ for each $P \in \mathcal{P}$;\footnote{Indeed, it will be assumed a symbol $\approx$ for the identity predicate such that $\varrho(\approx) =2$.} (ii)  a possible empty set $\mathcal{F}$ of  function  symbols, with the corresponding arity $\varrho(f) \geq 1$ for each $f \in \mathcal{F}$; (iii) a possible empty set of  individual constants $\mathcal{C}$. It will also assumed a fixed denumerable set $Var =\{x_1,x_2,\ldots\}$ of individual variables.\footnote{Most part of modal logic manuals --- for instance \cite{hug:cre:96}, \cite{fit:med:98} and \cite{gar:06} --- do not use function symbols among the symbols of their first-order modal signature (an exception is \cite[p. 241]{car:piz:08}). However, it is standard the use of function symbols in first-order signatures for classical logic, see for instance \cite[p. 49]{men:10} and \cite[p. 70]{ender:01}. It seems that this  absence in approaches to  first-oder modal logic is related to the problem of contingent identities, as we will discuss in Subsection~\ref{contingent}. We choose a wider approach by including function symbols, since in this new semantic the treatment of contingent identities is radically different comparing to Kripkean approach. As a result, some problems concerning references through the possible worlds will simply disappear.}

Give a first-order signature $\Theta$, the set $Ter(\Theta)$ of terms over $\Theta$ is defined recursively as follows: (i) a variable or a constant is a term; (ii) if $f$ is a $n$-ary function and $\tau_1, \ldots,  \tau_n$ are terms, then $f \tau_1 \ldots \tau_n$ is a term. The set $For(\Theta)$ of well-formed formulas (wffs) is also defined recursively as follows: (i) for each $n$-ary predicate $P$, 
if $\tau_1, \ldots, \tau_n$ are terms, then $P\tau_1\ldots \tau_n$ is a wff (called atomic); in particular,  if $\tau_1$ and $\tau_2$ are terms then $\approx\tau_1 \tau_2$, which will be written as $(\tau_1\approx \tau_2)$, is an atomic wff; (ii) if $\alpha$ is a wff and $x$ is a variable, then $(\neg \alpha)$, $(\Box \alpha)$ and $(\forall x \alpha)$ are also wffs; (iii) if $\alpha$ and $\beta$ are wffs, then $(\alpha \rightarrow \beta)$ is also a wff; (iv) nothing 
else is a wff. We will eliminate parenthesis when the readability is unambiguous.

The following abbreviations will be used in  {\bf Tm$^*$}:  $\alpha \vee \beta := \neg \alpha \to \beta$ (disjunction); $\alpha \wedge \beta := \neg(\alpha \to \neg\beta)$ (conjunction); $\Diamond \alpha := \neg\Box\neg\alpha$ (possibility); $\exists x \alpha :=\neg\forall x \neg\alpha$ (existential quantifier).
Let $\alpha$ be a wff of a first-order modal language $For(\Theta)$ over a first-order signature $\Theta$. The notion of free and bounded occurrences of a variable in a formula, as well as the notion of term free for a variable in a formula, closed term (that is, without variables) and closed formula (or sentence, that is, a formula without free occurrences of variables) are defined as usual (see, for instance, \cite{men:10}). The set of closed formulas and closed terms over $\Theta$ will be denoted by $Sen(\Theta)$ and $CTer(\Theta)$, respectively.  We write $\alpha[x/\tau]$ to denote the formula obtained from $\alpha$ by replacing simultaneously every free occurrence of the variable $x$ by the term $\tau$. If $\alpha$ is a formula and $y$ is a variable free for
the variable $x$ in $\alpha$,  $\alpha[x\wr y]$ denotes any
formula obtained from $\alpha$ by replacing some, but not
necessarily all (maybe none), free occurrences of $x$ by $y$.

\begin{definition} [da Costa] Let $\varphi$ and $\psi$ be formulas. If $\varphi$ can be obtained from $\psi$ by means of addition or deletion of void quantifiers,\footnote{That is, a quantifier $\forall x \alpha$ such that $x$ does not occur free in $\alpha$.} or by renaming bound variables (keeping the same free variables in the same places), we say that $\varphi$ and $\psi$ are {\em variant} of each other.
\end{definition}

\begin{definition} [Hilbert calculus for  {\bf Tm$^*$}] \label{axTm*}
The first-order modal logic {\bf Tm$^*$} is defined by a Hilbert calculus which consists of the following  axiom schemas and inference rules:\footnote{Recalling that $\Diamond \alpha$ is an abbreviation for $\neg\Box\neg\alpha$.}  \\

$\begin{array}{ll}
\axu & \alpha \rightarrow (\beta \rightarrow \alpha)\\[2mm]
\axd & (\alpha \rightarrow (\beta \rightarrow \gamma)) \rightarrow ((\alpha \rightarrow \beta) \rightarrow (\alpha \rightarrow \gamma))\\[2mm]
\axt & (\neg \beta \rightarrow \neg \alpha) \rightarrow ((\neg \beta \rightarrow \alpha) \rightarrow \beta)\\[2mm]
\axq & \forall x \alpha \rightarrow \alpha[x/\tau] \  \ \mbox{ if $\tau$ is free for $x$ in $\alpha$}\\[2mm]
\axc & \forall x (\alpha \rightarrow \beta) \rightarrow (\alpha \rightarrow \forall x\beta) \ \  \ \mbox{ if $\alpha$ contains no free  occurrences of $x$}\\[2mm]
\axs & \alpha\to\beta \ \ \mbox{ if $\alpha$ is a variant of $\beta$}\\[2mm]
\axst & \forall x(x \approx x)\\[2mm]
\axo & (x \approx y) \rightarrow (\alpha \rightarrow \alpha[x \wr y]) \ \   \ \mbox{ if $y$ is a variable free for $x$ in $\alpha$}\\[2mm]
\axNeq & (x \approx y) \rightarrow \Box (x \approx y)\\[2mm]
\axPeq & \neg(x \approx y) \rightarrow \Box \neg(x \approx y)\\[2mm]
\axK & \Box (\alpha \rightarrow \beta) \rightarrow (\Box \alpha \rightarrow \Box \beta)\\[2mm]
\axKu & \Box (\alpha \rightarrow \beta) \rightarrow (\Diamond  \alpha \rightarrow \Diamond\beta)\\[2mm]
\axKd & \Diamond(\alpha \rightarrow \beta) \rightarrow (\Box \alpha \rightarrow \Diamond \beta)\\[2mm]
\axMu & \Box \neg \alpha \rightarrow \Box(\alpha \rightarrow \beta)\\[2mm]
\axMd & \Box \beta \rightarrow \Box(\alpha \rightarrow \beta)\\[2mm]
\axMt & \Diamond \beta \to \Diamond (\alpha \rightarrow \beta) \\[2mm]
\axMc & \Diamond\neg \alpha \to \Diamond (\alpha \rightarrow \beta) \\[2mm]
\axT & \Box \alpha \rightarrow \alpha\\[2mm]
\axDNu & \Box \alpha \rightarrow \Box \neg \neg \alpha\\[2mm]
\axDNd & \Box \neg \neg \alpha \rightarrow \Box \alpha\\[2mm]
\axBF & \forall x \Box \alpha \rightarrow \Box \forall x \alpha\\[2mm]
\axCBF & \Box \forall x \alpha \rightarrow \forall x \Box \alpha\\[2mm]
\axNBF & \forall x \Diamond  \alpha \rightarrow \Diamond \forall x \alpha \\[2mm]
\axPBF &  \Diamond \forall x \alpha \rightarrow \forall x \Diamond  \alpha\\[4mm]
{\bf (MP)}: &  \beta \ \  \ \mbox{ follows from $\alpha$ and $\alpha \rightarrow \beta$}\\[2mm]
{\bf (Gen)}: &  \forall x \alpha \ \  \ \mbox{ follows from $\alpha$}
\end{array}
$
\end{definition}

The notion of derivation in {\bf Tm}$^*$ is defined as usual.  We will use the conventional notation $\Gamma \vdash_{{\bf Tm}^*} \alpha$ in order to express that there is a derivation in {\bf Tm}$^*$ of $\alpha$ from $\Gamma$.

As it could be expected, given that the rule of necessitation is not present in {\bf Tm}$^*$, this logic satisfies the restricted version of the Deduction metatheorem (DMT), as usually presented in first-order logics (see~\cite{men:10}):

\begin{theorem}[Deduction Metatheorem (DMT) for  {\bf Tm}$^*$]\label{DM}
Suppose that there exists in {\bf Tm}$^*$ a derivation of $\psi$ from $\Gamma \cup\{\varphi\}$, such that no application of the rule (Gen) has, as its quantified variable, a free variable of $\varphi$ (in particular, this holds when $\varphi$ is a sentence). Then $\Gamma \vdash_{{\bf Tm}^*} \varphi \to \psi$.\footnote{A more general version for quantified classical logic, which also holds for  {\bf Tm}$^*$, can be found in~\cite[p. 67]{men:10}.}
\end{theorem}

\begin{remark}  
It is well-known that normal modal logics admit two different notions of consequence relation with respect to Kripke semantics: the {\em local} one and the {\em global} one  (see, for instance, \cite[Defs. 1.35 and 1.37]{black-rij-ven:2002}). According to this, given a class $\mathbb{M}$ of Kripke models,
a formula $\varphi$ follows {\em locally} from a set $\Gamma$ of formulas if, for any $M \in \mathbb{M}$ and every world $w$ in $M$, $\varphi$ is true in $\langle M,w\rangle$ whenever every formula in $\Gamma$ is true in $\langle M,w\rangle$.  In turn,  $\varphi$ follows {\em globally} from $\Gamma$ in $\mathbb{M}$ if, for any $M \in \mathbb{M}$,  $\varphi$ is true in $\langle M,w\rangle$ for every $w$ whenever every formula in $\Gamma$ is true in $\langle M,w\rangle$ for every $w$.
From this, DMT holds in a propositional normal modal logic only for the local (semantical) consequence relation, and this is the perspective adopted with most normal modal logics, in which a modal logic can be studied simply in terms of validity of formulas. The same approach is assumed for first-order normal modal logics, ensuring the preservation of the deduction metatheorem. From the proof-theoretical perspective, the preservation of DMT even by considering the generalization inference rule {\bf (Gen)} and the necessitation rule forces to redefine the notion of derivation from premises in the corresponding Hilbert calculi for such modal logics. In this way, derivations are represented exclusively in terms of theoremhood (hence DMT holds by definition). Namely: $\Gamma \vdash \varphi$ iff either $\vdash \varphi$ or there exist $\gamma_1,\ldots,\gamma_n \in \Gamma$ such that $\vdash \gamma_1 \to(\gamma_2 \to(\ldots \to (\gamma_n \to \varphi)\ldots))$.

Concerning Ivlev-like non-normal modal first-order logics such as  {\bf Tm$^*$}, DMT holds exactly under the restrictions imposed in classical first-order logic w.r.t. the application of  {\bf (Gen)} to the premises. This is a consequence of assuming for  {\bf Tm$^*$} (as well as for the other first-order systems to be considered in this paper) the usual notion of derivation from premises in Hilbert calculi (which, in the case of first-order logics such as classical first-order logic, only works for global semantics).\footnote{Of course it would be possible to consider the notion of derivation from premises in the Hilbert calculi in which DMT holds by definition, as described above for standard modal logic.}

Another consequence of discarding the necessitation rule in the Ivlev-like systems is that the Replacement Metatheorem, which holds in every normal modal logic based on (propositional or first-order) classical logic, is no longer valid. That is, these logic are not self-extensional.\footnote{The Replacement Metatheorem says that if $\beta$ is a subformula of $\alpha$, $\alpha'$ is the result of replacing zero or more occurences of $\beta$ in $\alpha$ by a wff $\gamma$, then:  $\vdash \beta \leftrightarrow \gamma$ implies that $\vdash \alpha \leftrightarrow \alpha'$ (see, for instance~\cite[p. 72]{men:10}). In~\cite{omo:sku:16} it was shown that Replacement does not hold, in general, in Ivlev-like systems.}
\end{remark}

\subsection{Four-valued non-deterministic semantics for the logic {\bf Tm}$^*$} \label{structures}

In this section a suitable semantics of first-order structures for  {\bf Tm}$^*$ will be provided, formalizing the intuitions given in Section~\ref{intuition}. As we shall see in Remark~\ref{equiv-sem}, this semantics follows the lines of the non-deterministic first-order structures for paraconsistent logics based on Nmatrices given in~\cite{CFG19} which, by its turn, is adapted from the standard algebraic approach to first-order logic.

\begin{definition} \label{structure}
Let $\Theta$ be a first-order signature. A first-order structure for {\bf Tm}$^*$ is a pair $\mathfrak{A}= \langle U, \cdot^\mathfrak{A} \rangle$, such that $U$ is a non-empty set and $\cdot^\mathfrak{A}$ is an interpretation function for the symbols of $\Theta$ defined as follows:
\begin{itemize}
	\item For each $n$-ary predicate $P$, $P^\mathfrak{A} = (\mathfrak{a}^\mathfrak{A}(P),\mathfrak{c}^\mathfrak{A}(P))$ such that $\mathfrak{a}^\mathfrak{A}(P) \subseteq U^{n}$ and $\mathfrak{c}^\mathfrak{A}(P) \subseteq U^{n}$; it will be required that $\mathfrak{a}^\mathfrak{A}(\approx) = \{(a,a) \ : \ a \in U \}$ and $\mathfrak{c}^\mathfrak{A}(\approx) = \emptyset$; 
%it will be required that $\mathfrak{a}^\mathfrak{A}(\approx) \cup \mathfrak{c}^\mathfrak{A}(\approx) = \{(a,a) \ : \ a \in U \}$;
  \item For each individual constant $c$, $c^\mathfrak{A}$ is an element of $U$;
	\item For each $n$-ary function $f$, $f^\mathfrak{A}$ is a function from $U^n$ to $U$.
\end{itemize}
\end{definition}

\begin{remark} \label{equiv-sem}
Within the semantics for first-order languages defined over algebraic structures, a predicate symbol of arity $n$ is  interpreted by means of a function $I(P): U^n \to A$, where $U$ is the domain of the semantical structure and $A$ is the domain of a given algebra of truth-values. This generalizes the standard Tarskian structures in which $I(P)$ is a subset of $U^n$, which can be represented by its characteristic function $I(P): U^n \to \{0,1\}$. This approach was adapted in~\cite{CFG19} for non-deterministic first-order structures for paraconsistent logics based on Nmatrices over multialgebras called {\em swap structures}. In such framework, a predicate symbol $P$ is interpreted as a function $I(P): U^n \to B$ where $B$ is the domain of a given swap structure. It is easy to see that the notion of semantical structures for {\bf Tm}$^*$ given in Definition~\ref{structure} is equivalent to the above mentioned approach. Indeed, let $\mathfrak{A}= \langle U, \cdot^\mathfrak{A} \rangle$ be a  first-order structure for {\bf Tm}$^*$. For every $n$-ary predicate symbol $P$ let $P_\mathfrak{A}:U^n \to V_4$ be the function such that $P_\mathfrak{A}^{-1}(T^+)= \mathfrak{a}^\mathfrak{A}(P)\setminus \,\mathfrak{c}^\mathfrak{A}(P)$;  $P_\mathfrak{A}^{-1}(C^+)= \mathfrak{a}^\mathfrak{A}(P)\cap \,\mathfrak{c}^\mathfrak{A}(P)$;  $P_\mathfrak{A}^{-1}(C^-)= \mathfrak{c}^\mathfrak{A}(P)\setminus \,\mathfrak{a}^\mathfrak{A}(P)$; and  $P_\mathfrak{A}^{-1}(F^-)= U^n \setminus \big(\mathfrak{a}^\mathfrak{A}(P)\cup \,\mathfrak{c}^\mathfrak{A}(P)\big)$. This defines a first-order structure as in~\cite{CFG19}. Conversely, let $\mathfrak{A}$ be a first-order structure for  {\bf Tm}$^*$ as  in~\cite{CFG19}, that is: it is an structure as in Definition~\ref{structure}, with the only difference that the $n$-ary predicates are interpreted as functions $P_\mathfrak{A}:U^n \to V_4$. Now, define $P^\mathfrak{A} = (\mathfrak{a}^\mathfrak{A}(P),\mathfrak{c}^\mathfrak{A}(P))$ such that $\mathfrak{a}^\mathfrak{A}(P) = P_\mathfrak{A}^{-1}(T^+) \cup P_\mathfrak{A}^{-1}(C^+)$ and $\mathfrak{c}^\mathfrak{A}(P)=P_\mathfrak{A}^{-1}(C^+) \cup P_\mathfrak{A}^{-1}(C^-)$. This gives origin to a first-order structure for  {\bf Tm}$^*$ as in Definition~\ref{structure}. Therefore,  both semantical approaches are equivalent.
\end{remark}

\begin{definition} \label{assigments}
An {\em assignment} over a first-order structure  $\mathfrak{A}$ for  {\bf Tm}$^*$ is a function $s:Var \to U$. Given two assignments $s$ and $s'$ over $\mathfrak{A}$ and a variable $x$, we say that $s$ and $s'$ are {\em $x$-equivalent}, denoted by $s \sim_x s'$, if $s(y)=s'(y)$ for every $y \in Var$ such that $y\neq x$. The set of assignments  over $\mathfrak{A}$ will be denoted by $\mathsf{A}(\mathfrak{A})$. The set of assignments which are $x$-equivalent to $s$ will be denoted by $\mathsf{E}_x(s)$. If $a \in U$, the assignment $s' \in \mathsf{E}_x(s)$ such that $s'(x)=a$ will be denoted by $s^x_a$. 
\end{definition}

Given $\mathfrak{A}$ and $s$, the denotation $\termvalue{\tau}^\mathfrak{A}_s$ of a term $\tau$ in $(\mathfrak{A},s)$ is defined recursively as follows: (i) $\termvalue{x}^\mathfrak{A}_s = s(x)$ if $x$ is a variable; (ii) $\termvalue{c}^\mathfrak{A}_s = c^\mathfrak{A}$ if $c$ is a constant; and $\termvalue{f\tau_1\ldots\tau_n}^\mathfrak{A}_s=f^\mathfrak{A}(\termvalue{\tau_1}^\mathfrak{A}_s,\ldots,\termvalue{\tau_n}^\mathfrak{A}_s)$. Note that  $\termvalue{\tau}^\mathfrak{A}_s \in U$ for every term $\tau$.

\begin{definition} \label{Tm-sem}
Given $\mathfrak{A}$, a {\em valuation} over $\mathfrak{A}$ is a function $v:For(\Theta) \times \mathsf{A}(\mathfrak{A}) \to V_4$ defined recursively as follows:
\begin{enumerate}
  \item For atomic wffs of the form $P\tau_1\ldots \tau_n$,
	\begin{itemize}
		\item[-] \ $v(P\tau_1\ldots \tau_n, s) = T^+$ iff
				$(\termvalue{\tau_1}^\mathfrak{A}_s,\ldots, \termvalue{\tau_n}^\mathfrak{A}_s) \in \mathfrak{a}^\mathfrak{A}(P) \setminus \,\mathfrak{c}^\mathfrak{A}(P)$;
		\item[-] \ $v(P\tau_1\ldots \tau_n, s) = C^+$ iff
				$(\termvalue{\tau_1}^\mathfrak{A}_s,\ldots, \termvalue{\tau_n}^\mathfrak{A}_s) \in \mathfrak{a}^\mathfrak{A}(P) \cap \,\mathfrak{c}^\mathfrak{A}(P)$;
		\item[-] \  $v(P\tau_1\ldots \tau_n, s) = C^-$ iff
				$(\termvalue{\tau_1}^\mathfrak{A}_s,\ldots, \termvalue{\tau_n}^\mathfrak{A}_s) \in \mathfrak{c}^\mathfrak{A}(P) \setminus \,\mathfrak{a}^\mathfrak{A}(P)$;
		\item[-] \   $v(P\tau_1\ldots \tau_n, s) = F^-$ iff
				$(\termvalue{\tau_1}^\mathfrak{A}_s,\ldots, \termvalue{\tau_n}^\mathfrak{A}_s) \in U^n \setminus \big(\mathfrak{a}^\mathfrak{A}(P) \cup \,\mathfrak{c}^\mathfrak{A}(P)\big)$.
	\end{itemize}
	\item For atomic wffs of the form $(\tau_1\approx \tau_2)$,
	\begin{itemize}
		\item[-] \  $v((\tau_1\approx \tau_2),s) = T^+$ iff $\termvalue{\tau_1}^\mathfrak{A}_s = \termvalue{\tau_2}^\mathfrak{A}_s$;
		\item[-] \  $v((\tau_1\approx \tau_2),s) = F^-$ iff $\termvalue{\tau_1}^\mathfrak{A}_s \neq \termvalue{\tau_2}^\mathfrak{A}_s$;
	\end{itemize}	   
	   \item for wffs of the form $\# \alpha$ such that $\# \in \{\neg,\Box\}$, then $v(\# \alpha,s) \in \tilde{\#} v(\alpha,s)$, where $\tilde{\#}$ denotes the corresponding multioperation of $\mathcal{A}_{\bf Tm}$ associated to $\#$  as described in Subsection~\ref{Nmatrix}.

	   \item For wffs of the form $(\alpha \to \beta)$ then $v((\alpha \to \beta),s) \in v(\alpha,s) \,\tilde{\to}\, v(\alpha,s)$, where  $\tilde{\to}$ denotes the multioperation of $\mathcal{A}_{\bf Tm}$  as described in Subsection~\ref{Nmatrix}.
	   
	   \item For wffs of the form  $\forall x  \alpha$, let $X(\alpha,x,v,s) = \{v(\alpha,s') \ : \  s' \in \mathsf{E}_x(s)\}$ (recall Definition~\ref{assigments}). Then $v(\forall x \alpha,s) \in \tilde{\forall}_4^d\big(X(\alpha,x,v,s)\big)$, where  $\tilde{\forall}_4^d$ is defined as in Subsection~\ref{Nmatrix}.

	\item Let $\tau$ be a term free for a variable $z$ in formula $\varphi$, and let $b = \termvalue{\tau}^{\mathfrak{A}}_{s}$. Then $v(\varphi[z/\tau],s) = v(\varphi,s^z_b)$ (recalling the notation from Definition~\ref{assigments}).

\item If $\varphi$ and $\varphi'$ are variant, then $v(\varphi,s) = v(\varphi',s)$ for every $s$.

\item If $y$ is a variable free for $x$ in $\varphi$ then $v((x \approx y) \rightarrow (\varphi \rightarrow \varphi[x \wr y]),s) \in \{T^+,C^+\}$.
\end{enumerate}
\end{definition}

\begin{remark}
Recall the sentences (1)-(3) used in Section~\ref{intuition}  as motivating examples. Clause \emph{1.} in Definition~\ref{Tm-sem} intends to capture formally the informal considerations about the sentence (1) given above. In turn, clause \emph{2.} intends to formalize the considerations above about sentence (3).  Clause \emph{5.} is the formal counterpart of the considerations about the sentence (2).

The intuition about clauses \emph{3.} and \emph{4.} was already discussed in~\cite{con:cer:per:19}. The reader should ask if those non-deterministic operators, specially with respect to the implication operator, are rather arbitrary. In some sense, that is true and we agree with that criticism. We think, also, that even Ivlev would agree with that. Maybe this justifies that he proposed several different modal systems with respect to the implication operator. Anyway, the multioperator above for the implication is not completely arbitrary. In~\cite{con:cer:per:19} we presented some arguments in order to convince the reader that Ivlev's multioperator for {\bf Tm} implication is reasonably natural from an intuitive point of view. With respect to the negation, we don't have many options and we showed there that this operator is the most natural in this context. With respect to the semantical counterpart of the modal connective $\Box$, Ivlev proposed in~ \cite{ivl:88} not only the multioperator described in Subsection~\ref{Nmatrix}, but also many others. We choose this one because it does not collapse any iterations of the modal connective $\Box$. In Subsection~\ref{iterations} two cases of collapse of modal iterations will be discussed.

Finally, clauses \emph{6.} to \emph{8.} are necessary in order to deal within a non-deterministic framework in a coherent way, as it was already analyzed in~\cite{CFG19} for first-order paraconsistent logics. For instance, the substitution lemma (which is crucial in order to validate \axq) must be guaranteed by requiring to the valuations to  choose in a suitable way, and this is stated by clause~\emph{6.}. 
\end{remark}

%$\mathsf{A}(\mathfrak{A})$

\begin{definition}
Let $\mathfrak{A}$ be a first-order structure for {\bf Tm}$^*$, and let $v$ be a valuation over it. Then, $(\mathfrak{A},v)$ {\em satisfies} a wff $\alpha$ if $v(\alpha,s) \subseteq \{T^+,C^+\}$ for some  assignment $s$ over $\mathfrak{A}$. The formula $\alpha$ is {\em true} in 
$(\mathfrak{A},v)$ if it is satisfied by  every  assignment $s$. We say that $\alpha$  is valid in {\bf Tm}$^*$ if it is true in every $(\mathfrak{A},v)$ . Finally, a formula $\alpha$ is a {\em semantical consequence of a set $\Gamma$ of formulas w.r.t. first-order structures for {\bf Tm}$^*$}, denoted by $\Gamma \models_{{\bf Tm}^*} \alpha$, if, for every  structure $\mathfrak{A}$ and every valuation $v$, if every $\gamma \in \Gamma$ is true in $(\mathfrak{A},v)$  then $\alpha$ is true in $(\mathfrak{A},v)$.
\end{definition}

\begin{remark}
The notion of semantical entailment in {\bf Tm}$^*$ deserves some comments. Observe that the standard semantics for first-order classical logic w.r.t. Tarskian structures produces a unique denotation for formulas, given a structure and an assignment. The same holds for first-order structures defined over complete Boolean algebras, or over complete Heying algebras for first-order intuitionistic logic, among other examples of first-order algebraizable logics (see, for instance, the classical references~\cite{rasiowa,rasiowa-sikorski}). In the present non-deterministic semantics, it is possible to make several {\em choices} for the denotation of complex formulas in a given structure, which are performed by the valuations. This is why the basic semantical environment for {\bf Tm}$^*$ is given by pairs $(\mathfrak{A},v)$, and so the semantical consequence is given by taking assignments over such pairs, in order to interpret the free variables. The same situation happens with the non-deterministic semantics for paraconsistent logics presented in~\cite{CFG19}. 
\end{remark}

\subsection{Soundness}

In this section the soundness of {\bf Tm}$^*$ w.r.t. its non-deterministic semantics will be stated.

\begin{lemma} \label{free-var}
  Let $s$ and $s'$ be two assignments over $\mathfrak{A}$ which are $x$-equivalent. Then, $v(\alpha,s)=v(\alpha,s')$ for every formula $\alpha$ in which $x$ does not occur free.
\end{lemma}
\begin{proof}
It can be proved easily by induction on the complexity of $\alpha$.
\end{proof}

\begin{lemma} \label{A4}
Let $s$ be an assignment and $v$ a valuation over a structure $\mathfrak{A}$. Then:\\[1mm]
(1) $v(\forall x \alpha \to \alpha[x/\tau],s) \in \{T^+,C^+\}$ if $\tau$ is a term free for $x$ in $\alpha$.\\[1mm]
(2) $v(\forall x (\alpha \rightarrow \beta) \rightarrow (\alpha \rightarrow \forall x\beta),s)  \in \{T^+,C^+\}$ if $\alpha$ contains no free  occurrences of $x$.\\[1mm]
(3) $v(\forall x(x \approx x),s)=T^+$. \\[1mm]
(4) $v((x \approx y) \rightarrow \Box (x \approx y),s)  \in \{T^+,C^+\}$\\[1mm]
(5) $v(\neg(x \approx y) \rightarrow \Box \neg(x \approx y),s)  \in \{T^+,C^+\}$.\\[1mm]
(6) $v(\forall x \Box \alpha \rightarrow \Box \forall x \alpha,s)  \in \{T^+,C^+\}$.\\[1mm]
(7) $v(\Box \forall x \alpha \rightarrow \forall x \Box \alpha,s)  \in \{T^+,C^+\}$.\\[1mm]
(8) $v(\forall x \Diamond  \alpha \rightarrow \Diamond \forall x \alpha,s) \in \{T^+,C^+\}$. \\[1mm]
(9) $v(\Diamond \forall x \alpha \rightarrow \forall x \Diamond  \alpha.s)  \in \{T^+,C^+\}$.
\end{lemma}
\begin{proof}
(1) Suppose that $v(\forall x \alpha,s) \in \{T^+,C^+\}$. Then $\{v(\alpha,s') \ : \ s' \in \mathsf{E}_x(s)\} \subseteq \{T^+,C^+\}$, by Definition~\ref{Tm-sem}(5) and by the definition of $\tilde{\forall}_4^d$ given in Subsection~\ref{Nmatrix}. Let $b = \termvalue{\tau}^{\mathfrak{A}}_{s}$. Then $v(\alpha[x/\tau],s) = v(\alpha,s^x_b)$, by Definition~\ref{Tm-sem}(6). Since $s^x_b \in  \mathsf{E}_x(s)$ it follows that $v(\alpha,s^x_b) \in \{T^+,C^+\}$, by hypothesis. This means that $v(\alpha[x/\tau],s) \in \{T^+,C^+\}$  and so $v(\forall x \alpha \to \alpha[x/\tau],s) \in \{T^+,C^+\}$.\\[1mm]
(2) Assume that  $v(\forall x (\alpha \rightarrow \beta),s)  \in \{T^+,C^+\}$ and  $v(\alpha,s)  \in \{T^+,C^+\}$. As proved in (1), $\{v(\alpha \to \beta,s') \ : \ s' \in \mathsf{E}_x(s)\} \subseteq \{T^+,C^+\}$. Let $s' \in \mathsf{E}_x(s)$. Then $v(\alpha \to \beta,s') \in \{T^+,C^+\}$. But $v(\alpha \to \beta,s') \in v(\alpha,s') \,\tilde{\to}\, v(\beta,s')$, and $v(\alpha,s') = v(\alpha,s)$, by Lemma~\ref{free-var} (since $x$ does not occur free in $\alpha$). Then  $v(\alpha,s') \in \{T^+,C^+\}$ and so $v(\beta,s') \in \{T^+,C^+\}$. From this, $\{v(\beta,s')  \ : \ s' \in \mathsf{E}_x(s)\} \subseteq \{T^+,C^+\}$. This shows that $v(\forall x \beta,s)  \in \{T^+,C^+\}$, proving (2).\\[1mm]
(3) It is an immediate consequence of the definitions.\\[1mm]
(4) Suppose that $v((x \approx y),s)  \in \{T^+,C^+\}$. By Definition~\ref{Tm-sem}(2), $v((x \approx y),s) =T^+$ and so $v(\Box (x \approx y),s)  \in \{T^+,C^+\}$, by Definition~\ref{Tm-sem}(3) and by the definition of $\tilde{\Box}$.\\[1mm]
(5) Suppose that $v(\neg(x \approx y),s)  \in \{T^+,C^+\}$. By Definition~\ref{Tm-sem}(2), $v(\neg(x \approx y),s) =T^+$ and so (as in item~(4)) $v(\Box \neg(x \approx y),s)  \in \{T^+,C^+\}$.\\[1mm]
(6) Suppose that $v(\forall x \Box \alpha,s)  \in \{T^+,C^+\}$. Then $\{v(\Box \alpha,s')  \ : \ s' \in \mathsf{E}_x(s)\} \subseteq \{T^+,C^+\}$. By definition of $\tilde{\Box}$ (and by Definition~\ref{Tm-sem}(3)) it follows that $\{v(\alpha,s')  \ : \ s' \in \mathsf{E}_x(s)\} = \{T^+\}$. By Definition~\ref{Tm-sem}(5) and by the definition of $\tilde{\forall}_4^d$ it follows that $v(\forall x \alpha,s) = T^+$. From this $v(\Box \forall x \alpha,s)  \in \{T^+,C^+\}$, which proves (6).\\[1mm]
(7) Suppose that $v(\Box \forall x \alpha,s)  \in \{T^+,C^+\}$. Then $v(\forall x \alpha,s) = T^+$ and so $\{v(\alpha,s')  \ : \ s' \in \mathsf{E}_x(s)\} = \{T^+\}$. From this, $\{v(\Box \alpha,s')  \ : \ s' \in \mathsf{E}_x(s)\} \subseteq \{T^+,C^+\}$, therefore $v(\forall x \Box \alpha,s)  \in \{T^+,C^+\}$. This proves (7).\\[1mm]
(8) Suppose that $v(\forall x \Diamond \alpha,s)  \in \{T^+,C^+\}$. Then $\{v(\Diamond \alpha,s')  \ : \ s' \in \mathsf{E}_x(s)\} \subseteq \{T^+,C^+\}$. By definition of $\tilde{\Diamond}$ (and by Definition~\ref{Tm-sem}(3)) it follows that $F^- \not \in \{v(\alpha,s')  \ : \ s' \in \mathsf{E}_x(s)\}$. By Definition~\ref{Tm-sem}(5) and by the definition of $\tilde{\forall}_4^d$ it follows that $v(\forall x \alpha,s) \neq F^-$. From this $v(\Diamond \forall x \alpha,s)  \in \{T^+,C^+\}$, which proves (8).\\[1mm]
(9) Suppose that $v(\Diamond \forall x \alpha,s)  \in \{T^+,C^+\}$. Then $v(\forall x \alpha,s) \neq F^-$ and so $F^- \notin \{v(\alpha,s')  \ : \ s' \in \mathsf{E}_x(s)\}$. From this, $\{v(\Diamond \alpha,s')  \ : \ s' \in \mathsf{E}_x(s)\} \subseteq \{T^+,C^+\}$, therefore $v(\forall x \Diamond \alpha,s)  \in \{T^+,C^+\}$. This concludes the proof.
\end{proof}

\begin{lemma} \label{infrules}
Let  $\mathfrak{A}$ be a structure, and $v$ a valuation over it. Then:\\[1mm]
(1) If  $\alpha$ and $\alpha \rightarrow \beta$ are true in  $(\mathfrak{A},v)$, so is $\beta$.\\[1mm]
(2) If  $\alpha$ is true in  $(\mathfrak{A},v)$, so is $\forall x \alpha$.
\end{lemma}
\begin{proof}
It is an easy consequence of the definitions.
\end{proof}

\begin{theorem}[Soundness]
Let $\Gamma \cup\{\alpha\}$ be a set of formulas. Then: $\Gamma \vdash_{\bf Tm^*}\alpha$ implies that $\Gamma \models_{{\bf Tm}^*} \alpha$.
\end{theorem}
\begin{proof}
By definition of valuation, axioms \axs\ and \axo\ are true in every $(\mathfrak{A},v)$.
As it was proved in~\cite{con:cer:per:19}, all the schema axioms from {\bf Tm} are valid w.r.t. its four-valued  Nmatrix semantics, so they are also valid in ${\bf Tm^*}$. By Lemma~\ref{A4}, the rest of the schema axioms from ${\bf Tm^*}$ (not considered above) are also true in every $(\mathfrak{A},v)$. By Lemma~\ref{infrules}, the inference rules of ${\bf Tm^*}$ preserve trueness in a given structure and valuation. From this, and assuming that $\Gamma \vdash_{\bf Tm^*}\alpha$,  it is easy to prove, by induction on the length of a derivation in ${\bf Tm^*}$ of $\alpha$ from $\Gamma$, that $\Gamma \models_{{\bf Tm}^*} \alpha$.
\end{proof}	

\subsection{Completeness}

In this section, the completeness of ${\bf Tm^*}$ w.r.t. its structures will be obtained. In order to do this, the completeness proof for classical first-order logic with equality found in~\cite{chan:keis} will be adapted. This adaptation is similar to the one found in~\cite{CFG19} for first-order paraconsistent logics. Other details of the proof will be also based on the proof of completeness of the modal logic {\bf Tm} given in~\cite{con:cer:per:15} and~\cite{con:cer:per:17}.

\begin{definition}\label{Henkin}
Consider a theory $\Delta \subseteq For(\Theta)$ and a nonempty set $C$ of constants of the  signature $\Theta$. Then, $\Delta$ is called a
$C$-\emph{Henkin theory} in ${\bf Tm^*}$ if it satisfies the following: for every formula $\psi$   with (at most) a free variable $x$, there exists a constant $c$ in $C$ such that $\Delta \vdash_{\bf Tm^*} \psi[x/c] \to \forall x  \psi$.
\end{definition}

\begin{definition} \label{Tm-C} Let $\Theta_{C}$  be the signature obtained from $\Theta$ by adding a set $C$ of new individual constants. The consequence relation $\vdash^{C}_{\bf Tm^*}$ is the consequence relation of ${\bf Tm^*}$ over the signature $\Theta_{C}$.
\end{definition}

As it happens with first-order classical logic, the following result can be obtained for  ${\bf Tm^*}$:

\begin{theorem} \label{C-Hen}
Every  theory $\Delta \subseteq For(\Theta)$ in  ${\bf Tm^*}$ over a signature $\Theta$ can be conservatively extended to a $C$-Henkin theory $\Delta^H  \subseteq For(\Theta_C)$ in  ${\bf Tm^*}$ over a signature $\Theta_C$ as in Definition~\ref{Tm-C}. That is, $\Delta \subseteq \Delta^H$ and: $\Delta \vdash_{\bf Tm^*} \varphi$ iff $\Delta^H \vdash^C_{\bf Tm^*} \varphi$ for every $\varphi \in For(\Theta)$. In addition, if $\Delta^H \subseteq \overline{\Delta^H} \subseteq For(\Theta_C)$ then $\overline{\Delta^H}$ is also a $C$-Henkin theory.
\end{theorem}

Recall (see for instance Part~I of this paper~\cite{con:cer:per:19}) that, given a Tarskian and finitary logic  ${\bf L}=\langle  For,\vdash \rangle$ (where $For$ is the set of formulas of {\bf L}), and given a set of formulas $\Gamma \cup \{\varphi\} \subseteq For$, the set $\Gamma$ is  {\em maximal non-trivial with respect to $\varphi$} (or {\em $\varphi$-saturated}) in {\bf L} if the following holds:~(i)~$\Gamma \nvdash \varphi$, and~(ii)~$\Gamma,\psi \vdash\varphi$ for every $\psi \notin \Gamma$.  Observe that if $\Gamma$ is $\varphi$-saturated then, for every wff $\beta$: $\Gamma  \vdash \beta$ iff  $\beta \in \Gamma$. In addition, for any wff $\beta$: either $\beta \in \Gamma$ or $\neg \beta \in \Gamma$ (but not both simultaneously).

By adapting to  ${\bf Tm^*}$ a classical and general result by Lindenbaum and \L o\'s (see~\cite[Theorem~22.2]{woj:84} and~\cite[Theorem~4.8]{con:cer:per:19}) we obtain the following:

\begin{theorem} \label{saturated}
Let  $\Gamma \cup \{\varphi\} \subseteq For(\Theta)$ such that $\Gamma \nvdash_{{\bf Tm^*}} \varphi$. Then, there exists a set of formulas $\Delta \subseteq For(\Theta)$ which is $\varphi$-saturated in ${\bf Tm^*}$ and such that $\Gamma \subseteq \Delta$.
\end{theorem}

Now, a canonical structure will be defined.

\begin{proposition} \label{equrel}
Let $\Delta$ be a set of formulas over a signature $\Theta$ which is $\varphi$-saturated for a given formula $\varphi$ in ${\bf Tm^*}$. Suppose that, in addition, $\Delta$ is a $C$-Henkin theory for a set $C$ of constants in $\Theta$. Define over $C$ the following relation: $c \sim d$ iff $\Delta \vdash_{\bf Tm^*} (c \approx d)$. Then $\sim$ is an equivalence relation.
\end{proposition}
\begin{proof}
It is an immediate consequence of the properties of the identity predicate in  ${\bf Tm^*}$.
\end{proof}

For every $c \in C$ let $\tilde{c}=\{d \in C \ : \ c \sim d\}$ and $U_\Delta=\{\tilde{c} \ : \ c\in C\}$. 

\begin{definition} \label{canonstr}
Let $\Delta$, $\varphi$ and $C$ as in Proposition~\ref{equrel}. The canonical structure 
$\mathfrak{A}_\Delta = \langle U_\Delta, \cdot^{\mathfrak{A}_\Delta} \rangle$ is such that:

\begin{enumerate}[(a)]
	\item For each $n$-ary predicate $P$, $P^{\mathfrak{A}_\Delta} = (\mathfrak{a}^{\mathfrak{A}_\Delta}(P),\mathfrak{c}^{\mathfrak{A}_\Delta}(P))$ is such that:
	\begin{itemize}
	  \item[-] \ $(\tilde{c}_1, \ldots, \tilde{c}_n) \in \mathfrak{a}^{\mathfrak{A}_\Delta}(P) \textrm{ iff } 
	  P c_1 \ldots c_n \in \Delta$; 
		\item[-] \	$(\tilde{c}_1, \ldots, \tilde{c}_n) \in \mathfrak{c}^{\mathfrak{A}_\Delta}(P) 
			\textrm{ iff } \neg (\neg \Box  P c_1 \ldots c_n \rightarrow \Box \neg P c_1 \ldots
			c_n)  \in \Delta$.
	\end{itemize}
	\item For each $n$-ary function symbol, $f^{\mathfrak{A}_\Delta}$ is given by $f^{\mathfrak{A}_\Delta}(\tilde{c}_1, \ldots, \tilde{c}_n) := \tilde{c}$ where  $c \in C$ is such that $(fc_1 \ldots c_n \approx c) \in \Delta$.
	\item  For each individual constant $c$, $c^{\mathfrak{A}_\Delta} := \tilde{d}$, where  $d \in C$ is such that $(c \approx d) \in \Delta$.  
\end{enumerate}
\end{definition}

\begin{proposition} \label{can-wd}
The canonical structure $\mathfrak{A}_\Delta = \langle U_\Delta, \cdot^{\mathfrak{A}_\Delta} \rangle$ is well-defined.
\end{proposition}
\begin{proof}
The proof is an easy adaptation of the one for classical first-order logic, see for instance~\cite[Lemma~2.12]{chan:keis}. With respect to the definition of the functions $\mathfrak{a}^{\mathfrak{A}_\Delta}(P)$ and $\mathfrak{c}^{\mathfrak{A}_\Delta}(P)$, it is enough to observe that if $(c_i \approx d_i) \in \Delta$ for $1 \leq i \leq n$ then: $\varphi[x_1/c_1, \cdots, x_n/c_n] \in \Delta$ iff $\varphi[x_1/d_1, \cdots, x_n/d_n] \in \Delta$, for every formula $\varphi$ whose free variables occur in the list $x_1,\ldots,x_n$.
\end{proof}

\begin{definition} \label{can-fun}
Let $\mathfrak{A}_\Delta = \langle U_\Delta, \cdot^{\mathfrak{A}_\Delta} \rangle$ be the canonical structure for a set  of formulas $\Delta$ as in Definition~\ref{canonstr}. The {\em canonical valuation $v_\Delta$  over $\mathfrak{A}_\Delta$}  is defined as follows, for any assignment $s$ such that $s(x_i)=\tilde{c}_i$ for every $i \geq 1$:
	\begin{enumerate}[(a)]
	  \item $v_\Delta(\alpha,s)=T^+$ iff $\Box \alpha[x_1/c_1, \cdots, x_n/c_n] \in \Delta$; 
		\item $v_\Delta(\alpha,s)=C^+$ iff $\alpha[x_1/c_1, \cdots, x_n/c_n] \in \Delta$ and $\neg \Box \alpha[x_1/c_1, \cdots, x_n/c_n] \in  \Delta$; 
		\item $v_\Delta(\alpha,s)=C^-$ iff $\neg\alpha[x_1/c_1, \cdots, x_n/c_n] \in  \Delta$  and $\neg \Box \neg \alpha[x_1/c_1, \cdots, x_n/c_n]  \in  \Delta$; 
		\item $v_\Delta(\alpha,s)=F^-$ iff $\Box \neg \alpha[x_1/c_1, \cdots, x_n/c_n] \in  \Delta$,
	\end{enumerate}
where the free variables occurring in $\alpha$ belong to the list $x_1,\ldots,x_n$. 
\end{definition}

\begin{lemma} \label{can1}
 Let $\mathfrak{A}_\Delta = \langle U_\Delta, \cdot^{\mathfrak{A}_\Delta} \rangle$ be the canonical structure for a set  of formulas $\Delta$ as in Definition~\ref{canonstr}. Then,  $v_\Delta$ is a {\bf Tm$^*$}-valuation over $\mathfrak{A}_\Delta$. In addition,  for every formula $\beta \in For(\Theta)$  such that the variables occurring in it belong to the list $x_1,\ldots,x_n$, and for every assignment $s$ such that  $s(x_i)=\tilde{c}_i$ for every $i \geq 1$: $v_\Delta(\beta,s) \in \{T^+, C^+\}$ iff $\beta[x_1/c_1, \cdots, x_n/c_n] \in \Delta$.
\end{lemma}
\begin{proof}  First, observe that $v_\Delta$ is a well-defined function. It follows by the same arguments given in the proof of Proposition~\ref{can-wd}, by axiom \axT\ and by the properties of $\varphi$-saturated sets, namely: for every wff $\beta$, either $\beta \in \Delta$ or $\neg \beta \in \Delta$ (but not both simultaneously). The fact that, for every $\beta \in For(\Theta)$, $v_\Delta(\beta,s) \in \{T^+, C^+\}$ iff $\beta[x_1/c_1, \cdots, x_n/c_n] \in \Delta$ (for $s(x_i)=\tilde{c}_i$) is a consequence of the definition of $v_\Delta$ and the considerations above.
Thus, it suffices to prove that $v_\Delta$ satisfies the clauses 1-8 of Definition~\ref{Tm-sem}. The proof for clauses 1-5 of such definition will be done by induction on the complexity of the formula $\alpha$. Observe that, if $\tau$ is a term such that the variables occurring in it belong to the list $x_1,\ldots,x_n$ and $s(x_i)=\tilde{c}_i$ for every $i \geq 1$, then 
$$(\ast) \ \ \ \ \ \termvalue{\tau}^{\mathfrak{A}_\Delta}_s = \tilde{c} \ \ \mbox{ iff } \ \  (\tau[x_1/c_1, \cdots, x_n/c_n] \approx c) \in \Delta.$$

\begin{description}
  \item[Clause 1:] $\alpha$ is an atomic formula of the form $P\tau_1\ldots \tau_k$. Suppose that the variables occurring free in $\alpha$ belong to the list $x_1,\ldots,x_n$ and $s(x_i)=\tilde{c}_i$ for every $i \geq 1$. Let $\termvalue{\tau_i}^{\mathfrak{A}_\Delta}_s = \tilde{d_i}$ for $1 \leq i \leq k$.
	\begin{enumerate}[(a)]
	  \item Suppose that $v_\Delta(P\tau_1 \ldots \tau_k,s) = T^+$. Thus,
	  		 $(\Box P\tau_1 \ldots \tau_k)[x_1/c_1, \cdots, x_n/c_n] \in  \Delta$. By \axT, $(P\tau_1 \ldots \tau_k)[x_1/c_1, \cdots, x_n/c_n]  \in \Delta$. By $(\ast)$, $P d_1 \ldots d_k  \in \Delta$. By Definition~\ref{canonstr}(a),   		 
	  		  $(\tilde{d_1}, \ldots, \tilde{d_k}) \in \mathfrak{a}^{\mathfrak{A}_\Delta}(P)$ and so $(\termvalue{\tau_1}^{\mathfrak{A}_\Delta}_s, \ldots, \termvalue{\tau_k}^{\mathfrak{A}_\Delta}_s) \in \mathfrak{a}^{\mathfrak{A}_\Delta}(P)$.
			Besides,  $(\Box \alpha \rightarrow (\neg \Box \alpha \rightarrow \beta))[x_1/c_1, \cdots, x_n/c_n] \in \Delta$, since it is an instance of a theorem of classical logic. Then 
			$(\neg \Box P\tau_1 \ldots \tau_k \rightarrow \Box \neg P\tau_1 \ldots \tau_k)[x_1/c_1, \cdots, x_n/c_n] \in \Delta$, by (MP). 
			So,	$\neg(\neg \Box P\tau_1 \ldots \tau_k\rightarrow \Box \neg P\tau_1 \ldots \tau_k)[x_1/c_1, \cdots, x_n/c_n] \notin\Delta$, by $\Delta$-ma\-xi\-ma\-li\-ty. By reasoning as above, we conclude that  $(\termvalue{\tau_1}^{\mathfrak{A}_\Delta}_s, \ldots, \termvalue{\tau_k}^{\mathfrak{A}_\Delta}_s) \notin \mathfrak{c}^{\mathfrak{A}_\Delta}(P)$.
								
			Conversely, suppose that $(\termvalue{\tau_1}^{\mathfrak{A}_\Delta}_s, \ldots, \termvalue{\tau_k}^{\mathfrak{A}_\Delta}_s) \in \mathfrak{a}^{\mathfrak{A}_\Delta}(P) \setminus\, \mathfrak{c}^{\mathfrak{A}_\Delta}(P)$. By reasoning as above,  $(P\tau_1\ldots \tau_k)[x_1/c_1, \cdots, x_n/c_n] \in \Delta$. Since $(\termvalue{\tau_1}^{\mathfrak{A}_\Delta}_s, \ldots, \termvalue{\tau_k}^{\mathfrak{A}_\Delta}_s) \notin \mathfrak{c}^{\mathfrak{A}_\Delta}(P)$ then, by $\Delta$-maximality, $(\neg \Box P\tau_1 \ldots \tau_k \rightarrow \Box \neg P\tau_1 \ldots \tau_k)[x_1/c_1, \cdots, x_n/c_n] \in \Delta$. Suppose now that $(\neg \Box P\tau_1 \ldots \tau_k)[x_1/c_1, \cdots, x_n/c_n] \in \Delta$. From this we conclude, by (MP), that	$(\Box \neg P\tau_1 \ldots \tau_k)[x_1/c_1, \cdots, x_n/c_n] \in \Delta$. But then, as a consequence of axiom \axT\ we get that $(\neg P\tau_1 \ldots \tau_k)[x_1/c_1, \cdots, x_n/c_n] \in \Delta$ and so $\Delta$ would be inconsistent. 
			Thus, $(\neg \Box P\tau_1 \ldots \tau_k)[x_1/c_1, \cdots, x_n/c_n] \notin \Delta$	and so, by maximality,  $(\Box P\tau_1 \ldots \tau_k)[x_1/c_1, \cdots, x_n/c_n] \in \Delta$. Therefore,
			$v_\Delta(P\tau_1 \ldots \tau_k,s) = T^+$.
			
		\item Observe that $v_\Delta(P\tau_1 \ldots \tau_k,s) = C^+$ iff 
			$$(P\tau_1\ldots \tau_k)[x_1/c_1, \cdots, x_n/c_n] \in \Delta \ \mbox{  and } \ (\neg \Box P\tau_1 \ldots \tau_k)[x_1/c_1, \cdots, x_n/c_n] \in \Delta.$$
In this case $(\termvalue{\tau_1}^{\mathfrak{A}_\Delta}_s, \ldots, \termvalue{\tau_k}^{\mathfrak{A}_\Delta}_s) \in \mathfrak{a}^{\mathfrak{A}_\Delta}(P)$.
			Suppose now, by absurd, that 
			$$(\neg \Box P\tau_1 \ldots \tau_k \rightarrow \Box \neg P\tau_1 \ldots \tau_k)[x_1/c_1, \cdots, x_n/c_n] \in \Delta.$$ 
			By (MP), it would follows that
			$(\Box \neg P\tau_1 \ldots \tau_k)[x_1/c_1, \cdots, x_n/c_n] \in \Delta$. By \axT, we would have 
			$(\neg P\tau_1 \ldots \tau_k)[x_1/c_1, \cdots, x_n/c_n] \in \Delta$	and so $\Delta$ would be inconsistent. From that, we infer that
			$(\neg \Box P\tau_1 \ldots \tau_k \rightarrow \Box \neg P\tau_1 \ldots \tau_k)[x_1/c_1, \cdots, x_n/c_n] \notin\Delta$
			and, so, by  $\Delta$-ma\-xim\-ali\-ty, 
			$\neg(\neg \Box P\tau_1 \ldots \tau_k \rightarrow \Box \neg P\tau_1\ldots \tau_k)[x_1/c_1, \cdots, x_n/c_n] \in\Delta$.
			 Thus, $(\termvalue{\tau_1}^{\mathfrak{A}_\Delta}_s, \ldots, \termvalue{\tau_k}^{\mathfrak{A}_\Delta}_s) \in \mathfrak{c}^{\mathfrak{A}_\Delta}(P)$.
			 
			Conversely, suppose now that  $(\termvalue{\tau_1}^{\mathfrak{A}_\Delta}_s, \ldots, \termvalue{\tau_k}^{\mathfrak{A}_\Delta}_s) \in \mathfrak{a}^{\mathfrak{A}_\Delta}(P) \cap \, \mathfrak{c}^{\mathfrak{A}_\Delta}(P)$. By reasoning as in the previous cases, it follows that $(P\tau_1 \ldots \tau_k)[x_1/c_1, \cdots, x_n/c_n] \in \Delta$ and $\neg(\neg\Box P\tau_1 \ldots \tau_k\rightarrow \Box \neg P\tau_1 \ldots \tau_k)[x_1/c_1, \cdots, x_n/c_n] \in \Delta$.
			Now, suppose that $(\Box P\tau_1 \ldots \tau_k)[x_1/c_1, \cdots, x_n/c_n] \in \Delta$. But
			$$(\Box P\tau_1 \ldots \tau_k \rightarrow 
			(\neg \Box P\tau_1 \ldots \tau_k \rightarrow \Box \neg P\tau_1 \ldots \tau_k))[x_1/c_1, \cdots, x_n/c_n] \in \Delta,$$ 
			since it is an instance of a theorem of classical propositional logic.
			Thus, by (MP), we would have that
			$(\neg \Box P\tau_1 \ldots \tau_k \rightarrow \Box \neg P\tau_1 \ldots \tau_k)[x_1/c_1, \cdots, x_n/c_n] \in \Delta$ 
			and so $\Delta$ would be inconsistent.
			Hence, 
			$(\Box P\tau_1 \ldots \tau_k)[x_1/c_1, \cdots, x_n/c_n] \notin \Delta$, and so 
			$$(\neg \Box P\tau_1 \ldots \tau_k)[x_1/c_1, \cdots, x_n/c_n] \in \Delta,$$ 
			by $\Delta$-maximality. Therefore,  $v_\Delta(P\tau_1 \ldots \tau_k,s) = C^+$.

		\item Note that $v_\Delta(P\tau_1 \ldots \tau_k,s) = C^-$ \ iff 
		$$(\neg P\tau_1 \ldots \tau_k)[x_1/c_1, \cdots, x_n/c_n]  \in \Delta \ \mbox{  and } \ (\neg \Box \neg P\tau_1 \ldots \tau_k)[x_1/c_1, \cdots, x_n/c_n]  \in \Delta.$$ 
				So
			$(P\tau_1 \ldots \tau_k)[x_1/c_1, \cdots, x_n/c_n] \notin \Delta$, by $\Delta$-maximality. Hence, 
			$(\termvalue{\tau_1}^{\mathfrak{A}_\Delta}_s, \ldots, \termvalue{\tau_k}^{\mathfrak{A}_\Delta}_s) \notin 
			\mathfrak{a}^{\mathfrak{A}_\Delta}(P)$.
			Suppose now that
			$(\neg \Box P\tau_1 \ldots \tau_k  \rightarrow \Box \neg P\tau_1 \ldots \tau_k)[x_1/c_1, \cdots, x_n/c_n]\in \Delta$. 
			As $(\neg \delta \rightarrow \beta) \rightarrow (\neg \beta \rightarrow \delta)$ is a theorem of classical logic,
			we have that
			$$(\neg \Box \neg P\tau_1 \ldots \tau_k \rightarrow \Box P\tau_1 \ldots \tau_k)[x_1/c_1, \cdots, x_n/c_n] \in \Delta.$$ Hence, $(\Box P\tau_1 \ldots \tau_k)[x_1/c_1, \cdots, x_n/c_n] \in \Delta$. Thus,  $(P\tau_1 \ldots \tau_k)[x_1/c_1, \cdots, x_n/c_n] \in \Delta$ 
			by \axT, and so $\Delta$ would be inconsistent. Therefore,  
			$$(\neg \Box \neg P\tau_1 \ldots \tau_k  \rightarrow \Box \neg P\tau_1\ldots \tau_k)[x_1/c_1, \cdots, x_n/c_n]\notin\Delta.$$
			We conclude that 
			$\neg(\neg \Box P\tau_1 \ldots \tau_k  \rightarrow 
			\Box \neg P\tau_1 \ldots \tau_k)[x_1/c_1, \cdots, x_n/c_n] \in \Delta$, by $\Delta$-maximality, and so
			 $(\termvalue{\tau_1}^{\mathfrak{A}_\Delta}_s, \ldots, \termvalue{\tau_k}^{\mathfrak{A}_\Delta}_s) \in 
			\mathfrak{c}^{\mathfrak{A}_\Delta}(P)$.

			Conversely, suppose that $(\termvalue{\tau_1}^{\mathfrak{A}_\Delta}_s, \ldots, \termvalue{\tau_k}^{\mathfrak{A}_\Delta}_s) \in \mathfrak{c}^{\mathfrak{A}_\Delta}(P) \setminus\, \mathfrak{a}^{\mathfrak{A}_\Delta}(P)$. By reasoning as in the previous cases, it follows that			
			$(P\tau_1 \ldots \tau_k)[x_1/c_1, \cdots, x_n/c_n]\notin \Delta$ and $\neg(\neg \Box\neg P\tau_1\ldots\tau_k\rightarrow \Box\neg P\tau_1 \ldots \tau_k)[x_1/c_1, \cdots, x_n/c_n]\in \Delta$.
			By $\Delta$-maximality, 
			$$(\neg P\tau_1 \ldots \tau_k)[x_1/c_1, \cdots, x_n/c_n] \in \Delta.$$
			If $(\Box \neg P\tau_1 \ldots \tau_k)[x_1/c_1, \cdots, x_n/c_n] \in \Delta$, then we would have that 
			$(\neg \Box \neg P\tau_1 \ldots \tau_k  \rightarrow \Box \neg P\tau_1 \ldots \tau_k)[x_1/c_1, \cdots, x_n/c_n] \in \Delta$ and so 
			$\Delta$ woud be inconsistent. Hence, 
			$$(\Box \neg P\tau_1 \ldots \tau_k)[x_1/c_1, \cdots, x_n/c_n] \notin \Delta$$ 
			and so, by $\Delta$-maximality, $(\neg \Box \neg P\tau_1 \ldots \tau_k)[x_1/c_1, \cdots, x_n/c_n] \in \Delta$.
			Therefore, we conclude $v_\Delta(P\tau_1 \ldots \tau_k,s) = C^-$.

		\item $v_\Delta(P\tau_1 \ldots \tau_k,s) = F^-$ iff   
			$(\Box \neg P\tau_1 \ldots \tau_k)[x_1/c_1, \cdots, x_n/c_n]  \in \Delta$. 
			Thus, by \axT, $(\neg P\tau_1 \ldots \tau_k)[x_1/c_1, \cdots, x_n/c_n] \in \Delta$ and
			$(P\tau_1 \ldots \tau_k)[x_1/c_1, \cdots, x_n/c_n] \notin \Delta$, by $\Delta$-maximality. Therefore $(\termvalue{\tau_1}^{\mathfrak{A}_\Delta}_s, \ldots, \termvalue{\tau_k}^{\mathfrak{A}_\Delta}_s) \notin \mathfrak{a}^{\mathfrak{A}_\Delta}(P)$. 
			By \axu\ and (MP),
			$$(\neg \Box P\tau_1\ldots \tau_k \rightarrow \Box 
			\neg P\tau_1\ldots \tau_k)[x_1/c_1, \cdots, x_n/c_n] \in \Delta.$$ 
			Since $\Delta$	is consistent,
			$\neg(\neg \Box P\tau_1\ldots \tau_k \rightarrow \Box \neg P\tau_1\ldots\tau_k)[x_1/c_1, \cdots, x_n/c_n]\notin \Delta$.
			From this, it follows that $(\termvalue{\tau_1}^{\mathfrak{A}_\Delta}_s, \ldots, \termvalue{\tau_k}^{\mathfrak{A}_\Delta}_s) \notin \mathfrak{c}^{\mathfrak{A}_\Delta}(P)$.

			Suppose now that $(\termvalue{\tau_1}^{\mathfrak{A}_\Delta}_s, \ldots, \termvalue{\tau_k}^{\mathfrak{A}_\Delta}_s) \notin \mathfrak{a}^{\mathfrak{A}_\Delta}(P) \cup \,\mathfrak{c}^{\mathfrak{A}_\Delta}(P)$.  As above, it is concluded that
			$(P\tau_1 \ldots \tau_k)[x_1/c_1, \cdots, x_n/c_n] \notin \Delta$ and also
			$$\neg(\neg \Box P\tau_1\ldots \tau_k \rightarrow \Box \neg P\tau_1\ldots \tau_k)[x_1/c_1, \cdots, x_n/c_n] 
			\notin \Delta.$$ 
			Thus 
			$$(\neg P\tau_1 \ldots \tau_k)[x_1/c_1, \cdots, x_n/c_n] \in \Delta \ \mbox{ and }$$ 
			$$(\neg \Box P\tau_1\ldots \tau_k \rightarrow \Box \neg P\tau_1\ldots \tau_k)[x_1/c_1, \cdots, x_n/c_n] \in \Delta,$$ 
			by $\Delta$-maximality. 
			If 	$(\Box P\tau_1 \ldots \tau_k)[x_1/c_1, \cdots, x_n/c_n] \in \Delta$, then	$$(P\tau_1\ldots \tau_k)[x_1/c_1, \cdots, x_n/c_n] \in \Delta,$$ 
			by \axT, and so $\Delta$ would be inconsistent. Thus
			$(\Box P\tau_1 \ldots \tau_k)[x_1/c_1, \cdots, x_n/c_n] \notin \Delta$ and so, by $\Delta$-maximality, 
			$(\neg \Box P\tau_1 \ldots \tau_k)[x_1/c_1, \cdots, x_n/c_n] \in \Delta$. Hence, by (MP),
			$(\Box \neg P\tau_1 \ldots \tau_k)[x_1/c_1, \cdots, x_n/c_n] \in \Delta$. Therefore,	$v_\Delta(P\tau_1 \ldots \tau_k,s) = F^-$.

	\end{enumerate}
	
	\item[Clause 2:] $\alpha$ is an atomic formula of the form $(\tau_1 \approx \tau_2)$. As in {\bf Clause 1}, suppose that the variables occurring in $\tau_1$ and $\tau_2$ belong to the list $x_1,\ldots,x_n$ and $s(x_i)=\tilde{c}_i$ for every $i \geq 1$. Let $\termvalue{\tau_i}^{\mathfrak{A}_\Delta}_s = \tilde{d_i}$ for $1 \leq i \leq 2$.
	
	\begin{enumerate}[(a)]
		\item  Suppose that $v_\Delta((\tau_1 \approx \tau_2),s) = T^+$. Thus,  $(\Box (\tau_1 \approx \tau_2))[x_1/c_1, \cdots, x_n/c_n] \in \Delta$ and so, by \axT, $(\tau_1 \approx \tau_2)[x_1/c_1, \cdots, x_n/c_n] \in \Delta$. Then, $(d_1 \approx d_2) \in \Delta$ and so $\termvalue{\tau_1}^{\mathfrak{A}_\Delta}_s = \termvalue{\tau_2}^{\mathfrak{A}_\Delta}_s$.

		Conversely, suppose that  $\termvalue{\tau_1}^{\mathfrak{A}_\Delta}_s = \termvalue{\tau_2}^{\mathfrak{A}_\Delta}_s$. Then, $(d_1 \approx d_2) \in \Delta$ and so $(\tau_1 \approx \tau_2)[x_1/c_1, \cdots, x_n/c_n] \in \Delta$. Then, by \axNeq, 
		$(\Box(\tau_1 \approx \tau_2))[x_1/c_1, \cdots, x_n/c_n] \in \Delta$. Thus,	$v_\Delta((\tau_1 \approx \tau_2),s) = T^+$.
		
		 \item Suppose that $v_\Delta((\tau_1 \approx \tau_2),s)= F^-$.
		 Thus, $(\Box \neg (\tau_1 \approx \tau_2))[x_1/c_1, \cdots, x_n/c_n] \in \Delta$. By \axT, 	$(\neg (\tau_1 \approx \tau_2))[x_1/c_1, \cdots, x_n/c_n] \in \Delta$. 
			As $\Delta$ is nontrivial, 
			$$(\tau_1 \approx \tau_2)[x_1/c_1, \cdots, x_n/c_n] \notin \Delta.$$ 
			Hence, $(d_1 \approx d_2) \notin \Delta$
		and so	$\termvalue{\tau_1}^{\mathfrak{A}_\Delta}_s \neq \termvalue{\tau_2}^{\mathfrak{A}_\Delta}_s$.
		 
		 Suppose now that  $\termvalue{\tau_1}^{\mathfrak{A}_\Delta}_s \neq \termvalue{\tau_2}^{\mathfrak{A}_\Delta}_s$. Thus, $(\tau_1 \approx \tau_2)[x_1/c_1, \cdots, x_n/c_n] \notin \Delta$ and so $(\neg(\tau_1 \approx \tau_2))[x_1/c_1, \cdots, x_n/c_n] \in \Delta$, by $\Delta$-maximality.  Then, by \axPeq, $(\Box \neg (\tau_1 \approx \tau_2))[x_1/c_1, \cdots, x_n/c_n] \in \Delta$. In other words,	$v_\Delta((\tau_1 \approx \tau_2),s)= F^-$.
		 
	\end{enumerate}

		\item[Clauses 3-4:] $\alpha$ is of the form $\neg \beta$, $\Box \beta$ or $\beta \rightarrow \gamma$. The proof is analogous to the one for {\bf Tm} given in~\cite[Lemma~4]{con:cer:per:15} (see also~\cite[Lemma~4]{con:cer:per:17}).
	
	\item[Clause 5:] $\alpha$ is of the form $\forall x \beta$. As above, suppose that the variables occurring free in $\alpha$ belong to the list $x_1,\ldots,x_n$ and $s(x_i)=\tilde{c}_i$ for every $i \geq 1$. Assume, without loss of generality, that $x \neq x_i$ for $1 \leq i \leq n$.
	
	\begin{enumerate}[(a)]
	  \item Suppose that $v_\Delta(\forall x \beta,s)=T^+$. Then,	  $(\Box \forall x \beta)[x_1/c_1, \cdots, x_n/c_n]  \in \Delta$. By \axCBF\ and (MP) we have that $(\forall x \Box \beta)[x_1/c_1, \cdots, x_n/c_n] \in \Delta$. Thus, by \axq, we have that $(\Box \beta)[x_1/c_1, \cdots, x_n/c_n][x/c] \in \Delta$,	for every $c \in C$. This means that $v_\Delta(\beta,s') = T^+$ for every $s' \in \mathsf{E}_x(s)$. Then, $v_\Delta(\forall x \beta,s) \in \tilde{\forall}_4^d\big(X(\beta,x,v_\Delta,s)\big)$, where the notation is as in Definition~\ref{Tm-sem}(5).
		
		\item Suppose that $v_\Delta(\forall x \beta,s)=C^+$. Thus, 
		$(\forall x \beta)[x_1/c_1, \cdots, x_n/c_n]  \in \Delta$ and, in addition,  $(\neg \Box \forall x \beta)[x_1/c_1, \cdots, x_n/c_n]  \in \Delta$.  By \axq,  $\beta[x_1/c_1, \cdots, x_n/c_n][x/c] \in \Delta$ for every $c \in C$. This means that $v_\Delta(\beta,s') \subseteq \{T^+,C^+\}$ for every  $s' \in \mathsf{E}_x(s)$. If $v_\Delta(\beta,s') =T^+$ for every  $s' \in \mathsf{E}_x(s)$ then $(\Box \beta)[x_1/c_1, \cdots, x_n/c_n][x/c]  \in \Delta$, for every $c \in C$. Let $c \in C$ such that $\Delta \vdash_{\bf Tm^*} \psi[x/c] \to \forall x  \psi$, for $\psi=(\Box \beta)[x_1/c_1, \cdots, x_n/c_n]$. From this, $(\forall x \Box \beta)[x_1/c_1, \cdots, x_n/c_n] \in \Delta$. By \axBF,  $(\Box\forall x \beta)[x_1/c_1, \cdots, x_n/c_n] \in \Delta$, hence $\Delta$ would be trivial. This shows that $v_\Delta(\beta,s') =C^+$ for some  $s' \in \mathsf{E}_x(s)$. Then, $v_\Delta(\forall x \beta,s) \in \tilde{\forall}_4^d\big(X(\beta,x,v_\Delta,s)\big)$, where the notation is as in Definition~\ref{Tm-sem}(5).
			
		\item Suppose that $v_\Delta(\forall x \beta,s)=C^-$. Then,  $$(\neg \forall x \beta)[x_1/c_1, \cdots, x_n/c_n]  \in \Delta \ \mbox{ and }$$
		$$(\neg \Box \neg  \forall x \beta)[x_1/c_1, \cdots, x_n/c_n]  \in \Delta.$$ 
		By \axPBF, $(\forall x\neg \Box \neg  \beta)[x_1/c_1, \cdots, x_n/c_n]  \in \Delta$. Since $\Delta$ is $C$-Henkin,  
		$$(\neg \beta)[x_1/c_1, \cdots, x_n/c_n][x/c']  \in \Delta$$ for some $c' \in C$. On the other hand, $(\neg \Box \neg  \beta)[x_1/c_1, \cdots, x_n/c_n][x/c]  \in \Delta$ for every $c \in C$, by \axq. In particular, $(\neg \Box \neg  \beta)[x_1/c_1, \cdots, x_n/c_n][x/c']  \in \Delta$. This means that $v_\Delta(\beta,s') = C^-$ for some  $s' \in \mathsf{E}_x(s)$ (namely, for $s'=s^x_{a}$ with $a=\tilde{c'}$). If $v_\Delta(\beta,s'') = F^-$ for some  $s'' \in \mathsf{E}_x(s)$ then $(\Box\neg  \beta)[x_1/c_1, \cdots, x_n/c_n][x/c'']  \in \Delta$ for some $c'' \in C$, hence $\Delta$ would be trivial. Therefore $v_\Delta(\beta,s'') \neq F^-$ for every  $s'' \in \mathsf{E}_x(s)$ and so $v_\Delta(\forall x \beta,s) \in \tilde{\forall}_4^d\big(X(\beta,x,v_\Delta,s)\big)$, where the notation is as in Definition~\ref{Tm-sem}(5).
		
		\item Suppose that $v_\Delta(\forall x \beta,s)=F^-$. 
		Thus, 
		$$(\neg \forall x \beta)[x_1/c_1, \cdots, x_n/c_n]  \in \Delta \ \mbox{ and }$$
		$$(\Box \neg \forall x \beta)[x_1/c_1, \cdots, x_n/c_n]  \in \Delta.$$
		As $(\gamma \rightarrow \neg \delta) \rightarrow (\delta \rightarrow \neg \gamma)$ is a theorem of classical logic we have, by \axNBF, that $(\neg \forall x \neg \Box \neg \beta)[x_1/c_1, \cdots, x_n/c_n]  \in \Delta$. Since $\Delta$ is a $C$-Henkin theory and 
		$(\neg \gamma \rightarrow \delta) \rightarrow (\neg \delta \rightarrow \gamma)$ is a theorem of classical logic, it follows that 
		$$(\Box \neg \beta)[x_1/c_1, \cdots, x_n/c_n][x/c'] \in \Delta$$ 
		for some constant $c'$ in $C$. This means that $v_\Delta(\beta,s') = F^-$ for some  $s' \in \mathsf{E}_x(s)$ (namely, for $s'=s^x_{a}$ with $a=\tilde{c'}$).  Therefore  $v_\Delta(\forall x \beta,s) \in \tilde{\forall}_4^d\big(X(\beta,x,v_\Delta,s)\big)$, where the notation is as in Definition~\ref{Tm-sem}(5).

	\end{enumerate}
	
		\item[Clause 6:] (Substitution) Let $\tau$ be a term free for a variable $z$ in a formula $\psi$, and let $b=\tilde{c} = \termvalue{\tau}^{\mathfrak{A}_\Delta}_{s}$. By \axo\ (the Leibinz rule), $v_\Delta(\psi[z/\tau],s) = v_\Delta(\psi[z/c],s) = v_\Delta(\psi,s^z_b)$.

\item[Clause 7:] (Variant) Let $\psi$ and $\psi'$ two formulas such that their free variables  belong to the list $x_1,\ldots,x_n$,  and $s(x_i)=\tilde{c}_i$ for every $i \geq 1$.  Observe that, if $\psi$ and $\psi'$ are variant, so are $\psi[x_1/c_1, \cdots, x_n/c_n]$ and  $\psi'[x_1/c_1, \cdots, x_n/c_n]$, as well as the pairs of formulas $(\#\psi)[x_1/c_1, \cdots, x_n/c_n]$ and  $(\#\psi')[x_1/c_1, \cdots, x_n/c_n]$, for $\# \in \{\neg,\Box, \neg\Box,\Box\neg,\neg\Box\neg \}$. On the other hand, if $\psi$ and $\psi'$ are variant then $\psi \in \Delta$ iff $\psi' \in \Delta$, by \axs. From this, $v_\Delta(\psi,s)=v_\Delta(\psi',s)$ and so $v_\Delta$ satisfies clause~7 of Definition~\ref{Tm-sem}.

\item[Clause 8:] (Leibniz rule) Suppose that $y$ is a variable free for $x$ in a formula $\psi$. By \axo, $(x \approx y) \rightarrow (\psi \rightarrow \psi[x \wr y]) \in \Delta$. Therefore $v_\Delta((x \approx y) \rightarrow (\psi \rightarrow \psi[x \wr y]),s) \in \{T^+,C^+\}$.
\end{description}

This concludes the proof.
\end{proof}

\begin{theorem}[Completeness]
Let $\Gamma \cup\{\alpha\} \subseteq For(\Theta)$ be a set of  formulas. Then: $\Gamma \models_{\bf Tm^*}\varphi$ implies that $\Gamma \vdash_{{\bf Tm}^*} \varphi$.
\end{theorem}
\begin{proof}
Let  $\Gamma \cup \{\varphi\} \subseteq For(\Theta)$ such that $\Gamma \nvdash_{\bf Tm^*} \varphi$. Then, by Theorem~\ref{C-Hen}, there exists a $C$-Henkin theory $\Delta^{H}$ over $\Theta_{C}$ in ${\bf Tm^*}$ for a nonempty set $C$ of new individual constants such that $\Gamma \subseteq \Delta^{H}$ and, for every $\alpha \in For(\Theta)$:  $\Gamma \vdash_{\bf Tm^*} \alpha$ iff $\Delta^{H} \vdash^{C}_{\bf Tm^*} \alpha$. From this, $\Delta^{H} \nvdash^{C}_{\bf Tm^*} \varphi$. Then, by Theorem~\ref{saturated}, there exists a set of formulas $\overline{\Delta^{H}}$ in $\Theta_{C}$ extending $\Delta^{H}$ which is $\varphi$-saturated in ${\bf Tm^*}$. By Theorem~\ref{C-Hen}, $\overline{\Delta^{H}}$ is also a $C$-Henkin theory over $\Theta_{C}$ in $\bf Tm^*$.

Consider now the canonical structure $\mathfrak{A}_{\overline{\Delta^{H}}}$ for $\overline{\Delta^{H}}$ over signature  $\Theta_{C}$, as in Definition~\ref{canonstr}, and let $v_{\overline{\Delta^{H}}}$ be the canonical valuation over $\mathfrak{A}_{\overline{\Delta^{H}}}$  as in Definition~\ref{can-fun}. Let $s$ be an assignment over $\mathfrak{A}_{\overline{\Delta^{H}}}$, and let $s(x_i)=\tilde{c}_i$ for every $i \geq 1$. If $\gamma \in \Gamma$ such that the variables occurring in it belong to the list $x_1,\ldots,x_n$ then, by (Gen) and by \axq, $\gamma[x_1/c_1, \cdots, x_n/c_n] \in \Delta$. From this, $v_{\overline{\Delta^{H}}}(\gamma,s)  \in \{T^+, C^+\}$. This shows that every $\gamma \in \Gamma$ is true in $(\mathfrak{A}_{\overline{\Delta^{H}}},v_{\overline{\Delta^{H}}})$. On the other hand,  $\varphi \notin \overline{\Delta^{H}}$. If $\varphi$ is a closed formula then  $\varphi$ is not true in $(\mathfrak{A}_{\overline{\Delta^{H}}},v_{\overline{\Delta^{H}}})$, by the last assertion of Lemma~\ref{can1}. Otherwise, let $n= Max \, \{ i \ : \ x_i \ \mbox{ occurs free in $\varphi$}\}$. Let $\psi_0=\forall x_1 \cdots \forall x_{n-1}\varphi$. By \axq, $\forall x_n\psi_0 \not\in \overline{\Delta^{H}}$. Since $\overline{\Delta^{H}}$ is a $C$-Henkin theory, there exists $c_n \in C$ such that $\psi_0[x_n/c_n] \not \in \overline{\Delta^{H}}$. Let $\psi_1=\forall x_1 \cdots \forall x_{n-2}\varphi[x_n/c_n]$. By reasoning as above, there exists $c_{n-1} \in C$ such that $\psi_1[x_{n-1}/c_{n-1}] \not \in \overline{\Delta^{H}}$. Inductively, it can be proven that there are constants $c_1,\ldots,c_n \in C$ such that $\varphi[x_1/c_1, \cdots, x_n/c_n] \notin \Delta$. Let $s_0$ be an assignment such that $s_0(x_i)=\tilde{c}_i$ for every $1 \leq i \leq n$. Then, by the last assertion of Lemma~\ref{can1},  $v_{\overline{\Delta^{H}}}(\varphi,s_0)  \notin \{T^+, C^+\}$, showing that $\varphi$ is not true in $(\mathfrak{A}_{\overline{\Delta^{H}}},v_{\overline{\Delta^{H}}})$.
Now, let $\mathfrak{A}$  be the reduct of $\mathfrak{A}_{\overline{\Delta^{H}}}$  to signature $\Theta$, and let $v:For(\Theta) \times \mathsf{A}(\mathfrak{A}) \to V_4$ be the restriction of $v_{\overline{\Delta^{H}}}$ to $For(\Theta)$, that is: $v(\alpha,s)=v_{\overline{\Delta^{H}}}(\alpha,s)$ for every $\alpha \in For(\Theta)$ and every $s$ in $\mathsf{A}(\mathfrak{A})=\mathsf{A}(\mathfrak{A}_{\overline{\Delta^{H}}})$. Clearly, $v$ is a valuation over  $\mathfrak{A}$ and $s_0$ is an assignment over it such that $v(\varphi,s_0)\notin \{T^+, C^+\}$, hence $\varphi$ is not  true in $(\mathfrak{A},v)$ (as noted above, if $\varphi$ is a closed formula then $s_0$ can be arbitrary), while  every $\gamma \in \Gamma$ is true in $(\mathfrak{A},v)$. This shows that $\Gamma \not\models_{\bf Tm^*}\varphi$.
\end{proof}

\section{Philosophical questions concerning first-order modal logic} \label{filosofia}

In this section a brief philosophical discussion will be considered concerning first-order modal languages, analyzing the r\^ole that the new semantics proposed in the previous sections can play with respect to this debate. In particular, the problem of the contingent identities and the distinction between \emph{de re} and \emph{de dicto} modalities will be discussed. Of course we do not defend that the modal non-deterministic semantic presented here can solve those problems: we just want to invite the reader to rethink them from a different perspective. We are aware that the philosophical discussion about quantified modal logics is deep and rich: in this paper, which is mostly of technical character, we are just scratching the surface concerning these complex philosophical issues. A more detailed conceptual discussion is left for a future work.

\subsection{Contingent identities} \label{contingent}

We saw that, in standard Kripkean models, identities involving proper names are always necessary, since they are \emph{rigid designators}.
But Kripke thesis is stronger than that: he also defends that, differently from proper names, identities involving definite descriptions are contingent.\footnote{In \cite[p. 41]{kri:81}.} Consider, for instance, Kripke's example:

\begin{center}
(4) \ \ \  Richard Nixon is the winner of the presidential elections of the United States in 1968.
\end{center}

The sentence (4) is an identity involving a proper name and a definite description. Kripke argues that, although (4) is true, it might be false. Indeed, Richard Nixon could have lost the presidential elections. In contrast, Richard Nixon is, necessarily, Richard Nixon.

The relation between proper names and definite description is not as simple as Kripke argues. As showed by Ruth C. Barcan,  the following sentence is a theorem of any Kripke's first-order modal system:\footnote{It was first shown in \cite{bar:47} that it holds for Lewis system 
{\bf S2} and {\bf S4}, but is easy to check that it holds in any Kripkean normal modal system.}

$$(x \approx y) \to \Box (x \approx y)$$

Let $f_1$ be the function  for ``the winner of the presidential elections of $x$ in 1968''; let $c_1$ be the constant for ``Richard Nixon''; and let $c_2$ be the constant for ``United States''. Thus, (4) can be symbolized by:

$$\textrm{(4') } \ \ \ (c_1 \approx f_1 c_2)$$

Therefore, by the theorem above and (MP), we infer from (4') the following:

$$\Box (c_1 \approx f_1 c_2)$$

This force us to accept that, in Kripkean systems, not only equalities involving proper names, but also  those involving definite descriptions  are necessary too, against Kripke's own thesis.\footnote{See footnote \ref{converse_rigidity}. }  Indeed, the fact that $x$ and $y$ may be replaced by any term makes the distinction between rigid designator (that is, constants) and definite description (that is, functions) totally irrelevant.

The point here is that some kind of equalities seems to be contingent from an intuitive perspective, as pointed out by Hughes and Cresswell. That is the case of the sentence:

\begin{center}
(5) \ \ \  The person who lives next door is the mayor.
\end{center}

Surely, (5) is contingent, since it is clearly possible that the person who lives next door might not have been the mayor. 
However, according to the theorem above, if (5) is true, then it is necessarily true.

There are two ways, however, to overcome this difficulty. The first one is suggested by Kripke himself and it is inspired on Russellian distinction between wider and narrow scope.\footnote{In Kripke's words: ``someone might say that the man who taught Alexander might not have taught Alexander; though it could not have been true that: the man who taught Alexander didn't teach Alexander. This is Russell's distinction of scope.'' See \cite[p. 62]{kri:81}. With respect to definite description and Russell's distinction of scope, see \cite{rus:05}. Russell uses the expressions ``primary occurrence'' and ``secondary occurrence''.} First of all, we define the {\em  uniqueness quantification} as follows:

$$\exists ! x \alpha \equiv_{def} \exists x (\alpha \wedge \forall y (\alpha[x/y] \to (x \approx y)))$$

Let $P_1$ be the predicate ``$x$ is a person who lives next door'' and $P_2$ be ``$x$ is a mayor''. Thus, (5) is symbolized by:

$$\textrm{(5') } \ \ \  \exists !x P_1 x \wedge \exists! x P_2 x \wedge \forall x (P_1 x \to P_2 x)$$

In this interpretation,\footnote{See \cite[p. 318--322]{hug:cre:96}.} (5') can be clearly contingently true in Kripkean models. Consider, for instance, the following model: in  the world $w_1$ the individual $d_1$ is the only one who is in the extension of $P_1$ and the only one who is in the extesion of $P_2$. But in the world $w_2$ related to $w_1$, the individual $d_1$ is the only one who is in the extension of $P_1$, the individual $d_2$ is the only one who is in the extension of $P_2$ and $d_1 \neq d_2$. Clearly, the sentence (5') is true in $w_1$ but false in $w_2$. Therefore, (5') is contingently true in $w_1$.

Analogously, let $P_3$ be the predicate ``$x$ is  winner of the presidential election of the United States in 1968''. Thus, (4) can be formalized by:

$$\textrm{(4'') } \ \ \  \exists !x (P_3 x \wedge (x \approx c_1))$$

Consider, now, the following model: $c_1$ is the name of the individual $d_1$; in  the world $w_1$, the individual $d_1$ is the only one who is in the extension of $P_3$. But in the world $w_2$ related to $w_1$, the individual $d_2$ is the only one who is in the extension of $P_3$ and $d_1 \neq d_2$. Clearly, the sentence (4'') is true in $w_1$ but false in $w_2$. Therefore, (4'') is contingently true in $w_1$.

The second path is through the intuitive Carnapian notion of \emph{intensional objects}.\footnote{Carnap in \cite[p. 1--39]{car:47} did not use the expression ``intensional object''. He said that terms and predicates of the language (that he calls ``predicators'') have an intension and an extension. While ``rational animal'' and ``feather biped'' have the same extension, that is, hold for the same individuals, they don't have the same meaning and, thus, they are different with respect to the intension. In contrast, ``rational animal' and ``human'' have the same intension and the same extension too. The intension of an individual constant is called individual concept. Predicators and individual constants are extensionally equivalent if they are equivalent in some description state, while the are intensionally equivalent if they are equivalent in all description state. As pointed out by Carnap, the notion of intension tries to capture the Fregean notion of sense in \cite{frege:48}.} From a more contemporaneous perspective, we could understand  Carnapian proposal as follows.\footnote{We are following here \cite[p. 344--347]{hug:cre:96}.}  Let us consider constant domains semantics, for instance. Those models will be enriched by a set $I$ of allowable intensional objects, such that an intensional object $i$ is a function such that for each possible word $w$, $i(w)$ is a member of the domain, saying that the intensional object is in $w$. Thus, if ``the person who lives next door'' is the function $i$ and ``the mayor'' is the function $i'$, we can say that (5) is contingently true because $i(w) = i'(w)$, but it could be at least some word $w'$ related to $w$ such that $i(w') \neq i'(w')$. Analogously, if ``the winner of the presidential election of the United States in 1968'' is the function $i$ and $d_1$ is the individual refereed by the constant $c_1$, we can say that (4) is contingently true because $i(w) = d_1$, but it could be at least some word $w'$ related to $w$ such that $i(w') \neq d_1$.

In contrast, in Ivlev semantics we can deal with contingent identities without taking any of those paths. For that, let us consider the following axioms:\\
 
$\begin{array}{ll}
\axCequ  & \neg\Box (x \approx y)\\[2mm]
\axCeqd  & \neg\Box \neg (x \approx y)
\end{array}$

\

Let call ${\bf Tm^*_c}$ the system obtained replacing \axNeq\ and \axPeq\ with \axCequ\ and \axCeqd. It is clear that  axioms \axCequ\ and \axCeqd\ are considered for the purpose of capturing the notion of contingent identities. From a semantical point of view,
we change in Definition~\ref{Tm-sem} clause~\emph{2.} by the following:
	\begin{itemize}
		\item[-] \  $v((\tau_1 \approx \tau_2),s) = C^+$ iff $\termvalue{\tau_1}^{\mathfrak{A}}_s = \termvalue{\tau_2}^{\mathfrak{A}}_s$;
		\item[-] \   $v((\tau_1 \approx \tau_2),s) = C^-$ iff $\termvalue{\tau_1}^{\mathfrak{A}}_s \neq \termvalue{\tau_2}^{\mathfrak{A}}_s$.
	\end{itemize}

With few adjustments in  Lemma~\ref{can1}, it is easy to obtain completeness for ${\bf Tm^*_c}$. Thus, the proof for this case now runs as follows (recalling the notation stated in the proof of  Lemma~\ref{can1}):

\begin{enumerate}[(a)]
		\item  Suppose that $v_\Delta((\tau_1 \approx \tau_2),s) = C^+$. Thus,  $(\tau_1 \approx \tau_2)[x_1/c_1, \cdots, x_n/c_n] \in \Delta$ and so $(d_1 \approx d_2) \in \Delta$. From this, $\termvalue{\tau_1}^{\mathfrak{A}_\Delta}_s = \termvalue{\tau_2}^{\mathfrak{A}_\Delta}_s$.

		Conversely, suppose that  $\termvalue{\tau_1}^{\mathfrak{A}_\Delta}_s = \termvalue{\tau_2}^{\mathfrak{A}_\Delta}_s$. Then, $(d_1 \approx d_2) \in \Delta$ and so 
		$$(\tau_1 \approx \tau_2)[x_1/c_1, \cdots, x_n/c_n] \in \Delta.$$ 
		By \axCequ, 
		$(\neg\Box(\tau_1 \approx \tau_2))[x_1/c_1, \cdots, x_n/c_n] \in \Delta$. Thus,	$v_\Delta((\tau_1 \approx \tau_2),s) = C^+$.
		
		 \item Suppose that $v_\Delta((\tau_1 \approx \tau_2),s)= C^-$.
		 Thus, $(\neg (\tau_1 \approx \tau_2))[x_1/c_1, \cdots, x_n/c_n] \in \Delta$. As $\Delta$ is nontrivial, $(\tau_1 \approx \tau_2)[x_1/c_1, \cdots, x_n/c_n] \notin \Delta$. Hence, $(d_1 \approx d_2) \notin \Delta$
		and so	$\termvalue{\tau_1}^{\mathfrak{A}_\Delta}_s \neq \termvalue{\tau_2}^{\mathfrak{A}_\Delta}_s$.
		 
		 Suppose now that  $\termvalue{\tau_1}^{\mathfrak{A}_\Delta}_s \neq \termvalue{\tau_2}^{\mathfrak{A}_\Delta}_s$. Thus, $(\tau_1 \approx \tau_2)[x_1/c_1, \cdots, x_n/c_n] \notin \Delta$ and so $(\neg(\tau_1 \approx \tau_2))[x_1/c_1, \cdots, x_n/c_n] \in \Delta$, by $\Delta$-maximality.  By \axCeqd, 
		 $$(\neg\Box \neg (\tau_1 \approx \tau_2))[x_1/c_1, \cdots, x_n/c_n] \in \Delta.$$ 
		 In other words,	$v_\Delta((\tau_1 \approx \tau_2),s)= C^-$.
	\end{enumerate}

The considerations above are enough to realize that, in Ivlev semantics, the problem of contingent identities is independent of both the Russellian distinction between narrow/wider scope and the Carnapian distinction between intentional/extensional objects. From a philosophical point of view, this means that it is possible to support that there are contingent identities without being compromised with the  Russellian theory of definite description or with the Carnapian concept of intentional objects.

The reader may be asking why we should choose between a modal logic where the identities are always necessary and another one where all them are always contingent. In fact, it seems that identities like those of sentence (5) are contingent, while identities in  an arithmetical context, for example, are always necessary. Indeed,  such a radical choice is unnecessary. We can just use two different symbols, one  for each kind of identity relation. For instance, if $\approx$ is the identity symbol for necessary identities, the symbol $\approx_c$ could be used for contingent identities. In this case, we should add to {\bf Tm$^*$} the versions of \axst, \axo, 
\axCequ\ and \axCeqd\ using the new symbol $\approx_c$. Mathematical statement could be formalized with $\approx$, while factual or empirical identities (or, in Kripkean theory, identities involving some non-rigid designator) should be formalized using $\approx_c$.

The reader might argue that the option of defining two types of equality is also available in possible worlds semantics. That's true.
However, besides considering two different equalities $\approx$ and $\approx_c$, in Ivlev semantics we can define a third one with the following meaning: two individuals are equal if and only if  necessarily they are identical (in the first sense stated above). For this, we define:

$$\tau_1 \approxeq \tau_2 := \Box (\tau_1 \approx \tau_2).$$

Semantically, this produces the following:
	\begin{itemize}
		\item[-] \  $v((\tau_1 \approxeq \tau_2),s) \in \{T^+, C^+\}$ iff $\termvalue{\tau_1}^{\mathfrak{A}}_s = \termvalue{\tau_2}^{\mathfrak{A}}_s$;
		\item[-] \   $v((\tau_1 \approxeq \tau_2),s) \in \{C^-, F^-\}$ iff $\termvalue{\tau_1}^{\mathfrak{A}}_s \neq \termvalue{\tau_2}^{\mathfrak{A}}_s$.
	\end{itemize}

Observe that this non-deterministic possibility is not available in the possible worlds semantics, since Krikpe's models are deterministic. This turns evident the wide expressive power of the present non-deterministic semantical framework for first-order modal logics.

\subsection{Iterations between modal operators} \label{iterations}

Let's recall the following definition presented in subsection \ref{Nmatrix}:
$$\Diamond \alpha := \neg \Box \neg \alpha $$
There are two simple ways to obtain new first-order alethical systems that extend {\bf Tm}$^*$. First, we can add  axioms in order to iterate $\Box$ and $\Diamond$. This alternative, from a semantical point of view, means constraining the set of values of a non-deterministic matrix that interprets $\Box$.

We will consider, for instance, two ways  to iterate the modal operators by means of the following axioms:\\
 
$\begin{array}{ll}
\axfour  & \Box \alpha \rightarrow \Box \Box \alpha\\[2mm]
\axfive  & \Diamond \Box \alpha \rightarrow \Box \alpha
\end{array}$

\

In order to validate axioms \axfour\ and \axfive, the valuation functions of Definition~\ref{Tm-sem} should preserve, respectively the following multioperators:

\begin{displaymath}
	\begin{array}{c c c}
	\axfour & \axfive\\

		\begin{array}{|c|c|c|}
			\hline	x 	& \tilde{\Box} x & \tilde{\Diamond} x\\
			\hline 	T^+ 		& \{T^+\}  	& \{T^+\} \\
			\hline  C^+ 		& \{C^-,F^-\}  & \{T^+,C^+\} \\
			\hline  C^-			& \{C^-,F^-\} 	& \{T^+,C^+\} \\
			\hline 	F^- 		& \{C^-,F^-\} 	& \{C^-,F^-\} \\
			\hline
		\end{array}
		&
		\hspace*{5mm}\begin{array}{|c|c|c|}
			\hline x	& \tilde{\Box} x  & \tilde{\Diamond} x\\
			\hline  T^+		& \{T^+\} & \{T^+\}\\
			\hline 	C^+  	& \{F^-\} & \{T^+\}\\
			\hline 	C^- 	& \{F^-\} & \{T^+\}\\
			\hline 	F^- 	& \{F^-\} & \{F^-\}\\
			\hline
		\end{array}
	\end{array}
\end{displaymath}

We saw in the last subsection that there are two different ways of dealing with equalities. From that and from the two axioms presented above, we can define four new systems:
\begin{itemize}
  \item[-] \  ${\bf T4m^*} = {\bf Tm^*} \cup \{\axfour\}$
  \item[-] \  ${\bf T45m^*} = {\bf T4m^*} \cup \{\axfive\}$
   \item[-] \  ${\bf T4m^*_c} = {\bf Tm^*_c} \cup \{\axfour\}$
  \item[-] \  ${\bf T45m^*_c} = {\bf T4m^*_c} \cup \{\axfive\}$
\end{itemize}

Again, with few adjustments in  Lemma \ref{can1}, we can obtain the completeness results for any of those four new systems.\footnote{Based on \cite{con:cer:per:15} and \cite{con:cer:per:17}.} In fact, for each kind of equality, we can generate 14 different modal systems.\footnote{Keeping the (N)matrices of negation and implication fixed, as argued in \cite{con:cer:per:19}.}
The details on the completeness proof are purely technical and not so difficult to obtain.

\subsection{ \emph{De re} and \emph{de dicto} collapse} \label{de_re}

Kneale draws our attention to the fact that the formulas $\Diamond \forall x \alpha$ and $\exists x \Box  \alpha$ are stronger than $\forall x \Diamond  \alpha$ and $\Box \exists x \alpha$, respectively.\footnote{See \cite[p. 622--633]{kne:62}.} In fact, considering ${\bf S5^*}$ as the constant domain expansion of the  propositional modal system {\bf S5}, we have:

$$\Diamond \forall x \alpha \vDash_{\bf S5^*}\forall x \Diamond \alpha \textrm{ but } \forall x \Diamond  \alpha \nvDash_{\bf S5^*} \Diamond \forall x \alpha $$
$$\exists x \Box   \alpha\vDash_{\bf S5^*}\Box \exists x \alpha \textrm{ but } \Box \exists x \alpha \nvDash_{\bf S5^*} \exists x \Box \alpha \, . $$

The author defends that the difference between $\exists x \Box  \alpha$ and $\Box \exists x   \alpha$ is the formal counterpart of Abelardus' distinction between necessity \emph{de rebus} and necessity \emph{de sensu}. This distinction was found later on by Peter of Spain as the modal statements \emph{de re} and \emph{de dicto}.

According to Kneale, Abelard says that assertions of necessity \emph{de rebus} and \emph{de sensu} seem to entail each other. It is interesting to note that this philosophical position attributed to Abelard cannot be expressed in Kripke's framework. The reason is that,  in ${\bf S5^*}$,  none of such necessities  entails each other. But ${\bf S5^*}$ models are the strongest ones within this framework, in the sense that if two formulas are logically independent w.r.t. those models, so they are also  independent in any other Kripkean model.

The most astonishing result showed by Tich\'y in~\cite{tic:73} is that the distinction between \emph{de re} and \emph{de dicto} formulas is not eliminable, even by supplementing {\bf S5} with the axiom schema expressing the so-called {\em Principle of Predication} of von Wright, as cited below:

\begin{description}
	\item[(PP)] Properties divide into two types: those whose belonging to an object is always either necessary or impossible and those whose belonging to an object is always contingent.
\end{description}

This principle is symbolized  by  Tich\'y as follows (see~\cite[p. 391]{tic:73}):
$$(PP^*) \hspace{6mm} \forall x(\Box \alpha(x) \vee \Box \neg \alpha(x)) \vee \forall x(\Diamond \alpha(x) \wedge \Diamond \neg \alpha(x))$$

\noindent
for any wff $\alpha(x)$ in which $x$ is the unique variable (possibly) occurring free. Tich\'y proved that \emph{de re} formulas are not eliminable in {\bf S5$^*$} + $(PP^*)$. As a consequence of this, since {\bf S5$^*$} + $(PP^*)$ is stronger than any Kripkean first-order standard  semantics, Tichy's result holds, as a corollary, for any Kripkean system.

In contrast, we can easily check that \emph{de re} and \emph{de dicto} modalities entail each other in {\bf Tm$^ *$}. This is a direct consequence of \axNBF, \axPBF\ and the definition of $\exists$ and $\Diamond$. However, it could be possible to separate these notions if we consider the first definition of the universal quantifier we propose here,  as being a (non-deterministic) {\bf Tm$^ *$}-conjunction, namely the  multioperator $\tilde{\forall}_4$, instead of the deterministic operator $\tilde{\forall}_4^d$ adopted for {\bf Tm$^ *$} (recall Subsection~\ref{Nmatrix}). Indeed, consider an unary predicate symbol $P$ and let $\mathfrak{A}$ be a first-order structure for {\bf Tm$^ *$} such that
$\emptyset \neq \mathfrak{a}^\mathfrak{A}(P) \subset \mathfrak{c}^\mathfrak{A}(P) = U$. Then, for every $v$ and $s$, $\{v(P(x), s') \ : \ s' \in \mathsf{E}_x(s)\} = \{C^-,C^+\}$. This produces the following scenario:

$$
\begin{array}{|c|c|c|c|c|}
				\hline  X & \tilde{\Diamond}(X) &  \tilde{\forall}_4(\tilde{\Diamond}(X)) &  \tilde{\forall}_4(X) & \tilde{\Diamond}(\tilde{\forall}_4(X))\\
				\hline \{C^-, C^+\} & \{T^+, C^+\} & \{C^+\} & \{F^-, C^-\}& V_4 \\
				\hline
			\end{array}$$

\

\noindent Therefore, by taking $v(\Diamond\forall x P x,s)=C^-$, the formula $\forall x \Diamond Px \to \Diamond\forall x Px$ will be refuted. That is, \axNBF\ will be invalidated. Of course this is equivalent to invalidate the consequence relation $\Box\exists x Px \models \exists x \Box Px$. Indeed, in the same situation as above, by considering the  first definition of the existential quantifier as being a (non-deterministic) {\bf Tm$^ *$}-disjunction, namely the  multioperator $\tilde{\exists}_4$ instead of  $\tilde{\exists}_4^d$, we will have:

$$
\begin{array}{|c|c|c|c|c|}
				\hline  X & \tilde{\exists}_4(X) &  \tilde{\Box}(\tilde{\exists}_4(X)) &  \tilde{\Box}(X) & \tilde{\exists}_4(\tilde{\Box}(X))\\
				\hline \{C^-, C^+\} & \{T^+, C^+\}& V_4 & \{F^-, C^-\} & \{C^-\} \\
				\hline
			\end{array}$$

\

\noindent and so, by taking $v(\Box\exists x Px,s)=C^+$ the modal inference $\Box\exists x Px \models \exists x \Box Px$ will be invalidated. Observe that, in order to axiomatize $\tilde{\forall}_4$ in {\bf Tm$^ *$}, it is enough to delete axiom   \axNBF\  from the axiomatization of {\bf Tm$^ *$} given in Definition~\ref{axTm*}. It is also worth noting that the non-deterministic quantifier $\tilde{\forall}_4$ (as well as its dual $\tilde{\exists}_4$) still satisfies the Barcan formulas \axBF\ and \axCBF\ (of course, in the case of the existential quantifier, both formulas must be expressed in terms of $\Diamond$). Thus, the use of this quantifier just blocks   \axNBF, while preserving both Barcan formulas. In particular, \axPBF\ is validated by the non-deterministic quantifier $\tilde{\forall}_4$. Indeed, it is easy to check that the situation concerning  \axPBF\ is as follows:

$$
\begin{array}{|c|c|c|}
				\hline  X  & \tilde{\Diamond}(\tilde{\forall}_4(X)) & \tilde{\forall}_4(\tilde{\Diamond}(X))\\
				\hline F^- \notin X & V_4 \mbox{ or  } \{T^+, C^+\} & \{C^+\} \\
				\hline  F^- \in X & \{C^-, F^-\} & \{F^-\} \\
				\hline				
			\end{array}$$

\

\noindent The definition of an intuitive notion of universal quantifier in {\bf Tm$^ *$} which also rejects  \axPBF\  deserves future research.

Some philosophers may defend that the distinction between \emph{de re} and \emph{de dicto} modalities is very important from a metaphysical point of view. Others, like Abelardo, could claim that there is no such  distinction. The point here is that Kripke's relational semantics is not neutral with respect to that polemics: it forces us to go against Abelard, that is, it forces us to admit a metaphysical distinction that, according to Quine, throws us into {\em the jungle of Aristotelian essentialism}.\footnote{See \cite{qui:55}. It seems that {\bf Tm}$^*$ semantics is immune to Quine's three criticisms to modal logic. The first of these criticisms is that, in the propositional context, being necessary would be equivalent to being tautological. But this is not the case in Ivlev semantics, for a very simple reason: we can prove that $\alpha \vee \neg \alpha$ is always true, but if the value of $\alpha$ is $C^+$, we could choose a false value for $\Box (\alpha \vee \neg \alpha)$. The second criticism concerns,  roughly speaking, to the fact that identities are always necessary. As we have seen, this is not necessarily the case in our approach. The third criticism concerns the modal commitment of a kind of essentialism that, as argued above, does not apply to our systems.} In contrast, {\bf Tm$^ *$} offers a formal interpretation of ``necessary'' such that, in Quinean terms, is much more civilized.

\subsection{Barcan Formulas} \label{sectbarcan}

Finally, a famous problem of quantified modal logic will be briefly discussed here from the perspective of the non-deterministic semantics: the so-called {\em Barcan formulas}  \axBF\ and \axCBF, which were already mentioned above.

In standard Kripke semantics with constant domains, the formula \axBF\ is valid in those models in which the relational accessibility between possible world is, at least, reflexive and symmetrical,\footnote{See \cite[p.244-249]{hug:cre:96}.} while \axCBF\ is not valid in general. On the other hand, there are variable domains in which neither \axBF\ nor \axCBF\ are valid.\footnote{In varying domain models, \axBF\ is valid if and only if they are anti-monotonic, that is, if a world $w$ is related to $w'$, then the domain of $w'$ is a subset of the domain of $w$. In contrast, \axCBF\ is valid if and only the models are monotonic, that is, if $w$ is related to $w'$, then the domain of $w$ is a subset of the domains of $w'$. See \cite[p. 108-112]{fit:med:98}.}
Thus, Kripke models offer an apparatus sufficiently refined to invalidate \axCBF\ and, in some cases, even \axBF.

As analyzed above, both laws \axBF\ and \axCBF\ hold in any model of {\bf Tm$^*$} --- even when considering non-deterministic quantifers in order to block \axNBF. From this, the reader could ask if Ivlev's semantics is less expressive than Kripke models with variable domains, as long as the former collapses formulas that are distinguishable in Kripke's approach. It will be argued that this kind of conclusion is, at the very least, hasty, and that the situation is a little more complex than it seems.

Given a wff $\alpha$, consider the following notation:

$${\Box^n} \alpha \equiv_{def} \underbrace{\Box \Box \ldots \Box}_{ \textit{n times}}\alpha$$

Now, let us consider now the following {\em generalized Barcan} and {\em  generalized converse Barcan} formulas:\\

$
\begin{array}{ll}
\axBFn & \forall x \Box^n \alpha \rightarrow \Box^n \forall x \alpha\\[2mm]
\axCBFn & \Box^n \forall x \alpha \rightarrow \forall x \Box^n \alpha
\end{array}$

\

In any modal logic that validates the necessitation rule and axiom  \axK, if \axBF\  holds, then \axBFn\ will hold too. Those considerations also apply to \axCBF\  and \axCBFn.

That is not the case in Ivlev's semantics. Consider, for instance, a signature with an unary predicate symbol $P$. Let $\mathfrak{A}$ be a structure for {\bf Tm$^ *$} with domain $U$ such that  $\mathfrak{a}^\mathfrak{A}(P) = U$ and $\mathfrak{c}^\mathfrak{A}(P) = \emptyset$. Observe that $v(Px,s)=T^+$ (hence $v(\forall x Px,s) = T^+$) for every $v$ and $s$.
Now, let $v$ be a valuation over $\mathfrak{A}$ such that $v(\Box Px,s) = T^+$ and $v(\Box \forall x Px,s) = C^+$, for any $s$. From this, $\{v(\Box\Box Px,s')  \ : \ s' \in \mathsf{E}_x(s)\}$ is a set of designated values, for any $s$. Hence  $v(\forall x\Box\Box Px,s)$ is designated, for any $s$. On the other hand,  $v(\Box\Box \forall x Px,s)$ is not designated, for every $s$. Then, $(\mathfrak{A},v)$  does not validate \axBFd. This argument can be easily generalized to any \axBFn\ for $n > 2$.

Consider now a valuation $v'$ over the structure $\mathfrak{A}$ defined above such that, for any $s$, $v(\Box Px,s) = C^+$ and $v(\Box \forall x Px,s) = T^+$. Thus, $\{v(\Box\Box Px,s')  \ : \ s' \in \mathsf{E}_x(s)\}$ is a set of non-designated values, for any $s$, and so  $v(\forall x\Box\Box Px,s)$ is not designated, for any $s$. On the other hand,  $v(\Box\Box \forall x Px,s)$ is designated, for every $s$. This shows that $(\mathfrak{A},v)$  does not validate \axCBFd\ and, again, the argument can be easily generalized to any \axCBFn\ for $n > 2$.

From a philosophical point of view, the formula $\forall x \Box Px$ can be interpreted as saying that all the objects of the domain have essentially the propriety $P$. As we saw in Subsection~\ref{de_re}, this is a case of the \emph{de re} modality. The formula $\Box \forall x Px$, in turn, is a \emph{de dicto} modality.\footnote{In this sense, essentialism is ``the doctrine that (at least some) objects have (at least some) essential properties. This characterization is not universally accepted, but no characterization is; and at least this one has the virtue of being simple and straightforward'' (see \cite{sep-essential-accidental}). According to the author, the modal characterization of the essencialism is ``$P$  is an essential property of an object $o$ just in case it is necessary that $o$ has $P$, whereas $P$ is an accidental property of an object $o$ just in case $o$ has $P$ but it is possible that $o$ lacks $P$''. Formally, if $c$ is an individual constant representing $o$, we symbolize ``$P$ is an essential propriety of $o$'' by $\Box Pc$, whereas ``$P$ is an accidental property of $P$'' is symbolized by $Pc \wedge \Diamond \neg Pc$.} 

An essentialist does not want the collapse of modalities \emph{de re} and \emph{de dicto}. That is why the non-deterministic semantics proposed here could not be interesting for this philosophical perspective. In contrast, this Nmatrix semantics separates {\em levels} of essentialism: someone could hold that there is no first-level essentialism, while accepting the second-level one: the fact that all object have a property as essentially essential does not imply that it is necessarily necessary that all the objects have this property.  This new philosophical perspective could be called a \emph{multi-leveled essentialism}.

It is worth noting that, for Kripke, the relation between worlds must be, at least, reflexive.\footnote{Segerberg  defends that, against 
``followers of Kripke's terminology'', normal modal logics should cover those models whose accessibility relation is empty, see \cite[p. 12]{seg:71}.}
 In this case,  axiom \axT\ holds, hence the collapse of \axBF\ and \axBFn\ is unavoidable, and the same holds for \axCBF\ and \axCBFn. In turn, the system {\bf Tm}$^*$ is an interesting case in which, even with axiom \axT\ holding, there is no collapse at all, as we have seen above.

Finally, let us consider the {\em strict implication}, which is defined as usual as

$$\alpha \strictif \beta \equiv_{def} \Box(\alpha \to \beta)$$

Let us consider the following version of \axBF\ and \axCBF:\footnote{In fact, those are the original formulas considered by Ruth C. Barcan in \cite{bar:46}.}\\

$
\begin{array}{ll}
\axBFe & \forall x \Box \alpha \strictif \Box \forall x \alpha\\[2mm]
\axCBFe & \Box \forall x \alpha \strictif \forall x \Box \alpha
\end{array}$

\

In Kripke semantics, \axBF\ implies \axBFe, by the necessitation rule.  But that is not the case with Nmatrix semantic: just take an instance of \axBF\ as contingently true in some $(\mathfrak{A},v)$.\footnote{That is, consider a pair $(\mathfrak{A},v)$ and an instance $\gamma$ of \axBF\ such that, for every $s$, $v(\gamma,s)=C^+$. Then, $\Box\gamma$ is an instance of  \axBFe\ such that $v(\Box\gamma,s)$ is not designated, for evey $s$.} The argument is analogous with respect to \axCBF\ and \axCBFe.

Because of this, it is not possible to support an essentialist view with respect to strict implication but not with respect to material implication in the usual Kripkean context. We could call the stricter rules \axBFe\ and \axCBFe\  as \emph{strict essentialism}.  In this sense, the non-deterministic semantics of {\bf Tm$^*$} can model a version of essentialism, namely  \axBF\ and \axCBF, which is  weaker than strict essentialism. 

For these reasons, we believe that the  non-deterministic semantics for first-order modal logics proposed here can uncover subtleties that were previously invisible to us. This is a natural consequence of the fact that our eyes are too accustomed to facing these problems through the glasses of the possible worlds semantics.

\section{Considering other Ivlev-like modal systems} \label{other-sys}

The previous sections were devoted to analyze first-order extensions of the four-valued Ivlev system {\bf Tm}.
In previous papers (\cite{con:cer:per:15,con:cer:per:17,con:cer:per:19}) we investigated other Ivlev-like modal propositional systems, characterized by four-valued, six-valued and eight-valued Nmatrices.

Concerning four-valued systems, the extension   to first-order languages of some of them, such as {\bf T4m} and {\bf T45m}, were mentioned in Subsection~\ref{iterations}. The soundness and completeness results, as well as all topics discussed above about {\bf Tm$^*$}, can be easily adapted to the systems {\bf T4m$^*$} and {\bf T45m$^*$}.

Recall from Section~\ref{intuition} that the four-valued non-deterministic semantics for {\bf Tm} can be described in terms of the concepts of \emph{actually true} and \emph{contingently true}. From this, the predicate symbols in a structure $\mathfrak{A}$ for {\bf Tm$^*$} are interpreted in terms of the mappings $\mathfrak{a}^\mathfrak{A}$ and $\mathfrak{c}^\mathfrak{A}$ (which, as observed in  Remark~\ref{equiv-sem}, is equivalent to consider a function $P_\mathfrak{A}:U^n \to V_4$, for every $n$-ary predicate symbol $P$). However, as discussed in Part~I of this paper (see~\cite[Section~1]{con:cer:per:19}), this is a particular case of a more general situation involving eight truth-values, each of one explained in terms of the modal concepts of \emph{necessarily true}, \emph{possibly true} and \emph{actually true}:

\begin{itemize}
  \item[] \ $T^+$: necessarily, possibly and actually true;
  \item[] \  $C^+$: contingently and actually true;
  \item[] \  $F^+$: impossible, possibly false but actually true;
  \item[] \  $I^+$: necessary true, impossible and actually true;
  \item[] \  $T^-$: necessarily and possibly true but actually false;
  \item[] \  $C^-$: contingently and actually false;
  \item[] \  $F^-$: impossible, possibly false and actually false;
  \item[] \  $I^-$: necessarily true, impossible and actually false. 
\end{itemize}

\noindent Under this perspective,  $+ = \{T^+,C^+,F^+,I^+ \}$ represents being {\em actually true} (the designated truth-values), while  being {\em actually false} is given by  $- = \{T^-,C^-,F^-,I^- \}$ (the undesignated truth-values). This produces an eight-valued non-normal version of {\bf K} that we called {\bf Km}. However, as observed in Part~I, the values $I^+$ and $I^-$ represent  very artificial situations. This lead us to consider six-valued structures, apt to deal with a deontic non-normal logic called {\bf Dm}, as well as some other axiomatic extensions of it.
The interpretation of a predicate symbol $P$ in eight-valued structures for {\bf Km} and in six-valued structures for {\bf Dm}, respectively,  are represented as follows:

\begin{center}
\begin{tikzpicture}[thick] 
\draw (2.7,-2.54) rectangle (-1.5,1.5) node[below right] {}; 
\draw (0,0) circle (1) node[above,shift={(0,1)}] {$\mathfrak{a}^\mathfrak{A}$};
\draw (1.2,0) circle (1) node[above,shift={(0,1)}] {$\mathfrak{n}^\mathfrak{A}$};
\draw (.6,-1.04) circle (1) node[shift={(1.1,-.6)}] {$\mathfrak{p}^\mathfrak{A}$}; 
	\node at (.6,-.4) {$T^+$}; 
	\node at (1.2,-.7) {$T^-$}; 
	\node at (0,-.7) {$C^+$}; 
	\node at (1.4,.2) {$I^-$}; 
	\node at (.6,.3) {$I^+$}; 
	\node at (-.2,.2) {$F^+$}; 
	\node at (.5,-1.4) {$C^-$}; 
	\node at (-1,-2) {$F^-$}; 
\end{tikzpicture}
\hspace*{1cm}\begin{tikzpicture}[thick] 
\draw (2.6,-1.8) rectangle (-2,1.8) node[below right] {}; 
\draw (-0.7,0) circle (1) node[above,shift={(0,1)}] {$\mathfrak{a}^\mathfrak{A}$};
\draw (0.6,0) ellipse (1cm and 0.5cm) node[above,shift={(0,0.4)}] {$\mathfrak{n}^\mathfrak{A}$};
\draw (0.6,0) ellipse (1.7cm and 1cm) node[shift={(1.5,0.8)}] {$\mathfrak{p}^\mathfrak{A}$}; 
	\node at (-1.4,0) {$F^+$}; 
	\node at (-0.7,0) {$C^+$}; 
	\node at (0,0) {$T^+$}; 
	\node at (0,0) {$T^+$}; 
	\node at (0.8,0) {$T^-$}; 
	\node at (2,0) {$C^-$};	
	\node at (-1.5,-1.3) {$F^-$}; 
\end{tikzpicture}
\end{center}

Recall from~\cite[Section~4]{con:cer:per:19} that the multioperations for the Nmatrix characterizing the logic {\bf Km} are as follows, where $\Diamond \alpha := \neg \Box \neg \alpha$, $\alpha \vee \beta := \neg \alpha \to \beta$ and $\alpha \wedge \beta := \neg(\alpha \to \neg\beta)$ and $+$ is the set of designated values:

\begin{displaymath}
\begin{array}{|l|c|c|c|}
   	   \hline  	   & \tilde{\neg}		&\tilde{\Box} 	& \tilde{\Diamond} \\
\hline \hline  T^+ & \{F^-\} 	& +		& + 	\\
 		\hline C^+ & \{C^-\} 	& -		& +  \\
		\hline F^+ & \{T^-\} 	& -		& - \\		
		\hline I^+ & \{I^-\} 	& +     & -	\\		
		\hline T^- & \{F^+\}	& + 	& + \\
 		\hline C^- & \{C^+\} 	& -		& + \\
 		\hline F^- & \{T^+\} 	& - 	& - \\
		\hline I^- & \{I^+\} 	& +     & -	\\		
	 	\hline	
\end{array}
\end{displaymath}

\

\begin{displaymath}
\begin{array}{|c|c|c|c|c|c|c|c|c|}
      \hline \tilde{\to} 	& T^+ 	& C^+ 		 	 & F^+ 	   & I^+ 	  & T^- 	& C^- 				& F^- 	  & I^-\\
\hline \hline     		T^+ & \{T^+\} & \{C^+\} 		 & \{F^+\} & \{I^+\} & \{T^-\} & \{C^-\} 			& \{F^-\} & \{I^-\}\\
       \hline     		C^+ & \{T^+\} & \{T^+,C^+\} 	 & \{C^+\} & \{I^+\} & \{T^-\} & \{T^-,C^-\} 		& \{C^-\} & \{I^-\}\\
       \hline     		F^+ & \{T^+\} & \{T^+\} 		 & \{T^+\} & \{I^+\} & \{T^-\} & \{T^-\} 			& \{T^-\} & \{I^-\}\\
		\hline			I^+ & \{I^+\} & \{I^+\} 		 & \{I^+\} & \{I^+\} & \{I^-\} & \{I^-\} 			& \{I^-\} & \{I^-\}\\    
       \hline     		T^- & \{T^+\} & \{C^+\} 		 & \{F^+\} & \{I^+\} & \{T^+\} & \{C^+\} 			& \{F^+\} & \{I^+\}\\
       \hline     		C^- & \{T^+\} & \{T^+,C^+\} 	 & \{C^+\} & \{I^+\} & \{T^+\} & \{T^+,C^+\} 		& \{C^+\} & \{I^+\}\\
       \hline     		F^- & \{T^+\} & \{T^+\} 		 & \{T^+\} & \{I^+\} & \{T^+\} & \{T^+\} 			& \{T^+\} & \{I^+\}\\
       \hline			I^- & \{I^+\} & \{I^+\} 		 & \{I^+\} & \{I^+\} & \{I^+\} & \{I^+\} 			& \{I^+\} & \{I^+\}\\    
       \hline
\end{array}
\end{displaymath}

\begin{displaymath}
\begin{array}{|c|c|c|c|c|c|c|c|c|}
      \hline \tilde{\vee} 	& T^+ 	& C^+ 		 	 & F^+ 	   & I^+ 	  & T^- 	& C^- 				& F^- 	  & I^-\\
\hline \hline
T^+ & \{T^+\} & \{T^+\} 		 & \{T^+\} & \{I^+\} & \{T^+\} & \{T^+\} 			& \{T^+\} & \{I^+\}\\
       \hline
C^+ & \{T^+\} & \{T^+,C^+\} 	 & \{C^+\} & \{I^+\} & \{T^+\} & \{T^+,C^+\} 		& \{C^+\} & \{I^+\}\\
       \hline
F^+ & \{T^+\} & \{C^+\} 		 & \{F^+\} & \{I^+\} & \{T^+\} & \{C^+\} 			& \{F^+\} & \{I^+\}\\
		\hline
I^+ & \{I^+\} & \{I^+\} 		 & \{I^+\} & \{I^+\} & \{I^+\} & \{I^+\} 			& \{I^+\} & \{I^+\}\\    
       \hline
T^- & \{T^+\} & \{T^+\} 		 & \{T^+\} & \{I^+\} & \{T^-\} & \{T^-\} 			& \{T^-\} & \{I^-\}\\
       \hline
C^- & \{T^+\} & \{T^+,C^+\} 	 & \{C^+\} & \{I^+\} & \{T^-\} & \{T^-,C^-\} 		& \{C^-\} & \{I^-\}\\
       \hline
F^- & \{T^+\} & \{C^+\} 		 & \{F^+\} & \{I^+\} & \{T^-\} & \{C^-\} 			& \{F^-\} & \{I^-\}\\
       \hline
I^- & \{I^+\} & \{I^+\} 		 & \{I^+\} & \{I^+\} & \{I^-\} & \{I^-\} 			& \{I^-\} & \{I^-\}\\    
       \hline
\end{array}
\end{displaymath}

\begin{displaymath}
\begin{array}{|c|c|c|c|c|c|c|c|c|}
      \hline \tilde{\wedge} 	& T^+ 	& C^+ 		 	 & F^+ 	   & I^+ 	  & T^- 	& C^- 				& F^- 	  & I^-\\
\hline \hline
T^+ & \{T^+\} & \{C^+\} 		 & \{F^+\} & \{I^+\} & \{T^-\} & \{C^-\} 			& \{F^-\} & \{I^-\}\\
       \hline
C^+ & \{C^+\} & \{F^+,C^+\} 	 & \{F^+\} & \{I^+\} & \{C^-\} & \{F^-,C^-\} 		& \{F^-\} & \{I^-\}\\
       \hline
F^+ & \{F^+\} & \{F^+\} 		 & \{F^+\} & \{I^+\} & \{F^-\} & \{F^-\} 			& \{F^-\} & \{I^-\}\\
		\hline
I^+ & \{I^+\} & \{I^+\} 		 & \{I^+\} & \{I^+\} & \{I^-\} & \{I^-\} 			& \{I^-\} & \{I^-\}\\    
       \hline
T^- & \{T^-\} & \{C^-\} 		 & \{F^-\} & \{I^-\} & \{T^-\} & \{C^-\} 			& \{F^-\} & \{I^-\}\\
       \hline
C^- & \{C^-\} & \{F^-,C^-\} 	 & \{F^-\} & \{I^-\} & \{C^-\} & \{F^-,C^-\} 		& \{F^-\} & \{I^-\}\\
       \hline
F^- & \{F^-\} & \{F^-\} 		 & \{F^-\} & \{I^-\} & \{F^-\} & \{F^-\} 			& \{F^-\} & \{I^-\}\\
       \hline
I^- & \{I^-\} & \{I^-\} 		 & \{I^-\} & \{I^-\} & \{I^-\} & \{I^-\} 			& \{I^-\} & \{I^-\}\\    
       \hline
\end{array}
\end{displaymath}

\

The Nmatrix for {\bf Dm} is obtained from this by eliminating the values $I^+$ and $I^-$, while the Nmatrix for {\bf Tm} is obtained by additionally removing the values  $F^+$ and $T^-$. Indeed, and as shown in~\cite[Proposition~5.5]{con:cer:per:19}, $\mathcal{A}_{\bf Tm} \subseteq_{sm} \mathcal{A}_{\bf Dm} \subseteq_{sm}\mathcal{A}_{\bf Km}$. Here, $\mathcal{A}_{\bf Dm}$ and $\mathcal{A}_{\bf Km}$ denote respectively the multialgebras for {\bf Dm} and {\bf Km}, and $\subseteq_{sm}$ is the `submultialgebra' relation.\footnote{By comparing the figure for $\mathfrak{a}^\mathfrak{A}(P)$ and $\mathfrak{c}^\mathfrak{A}(P)$ displayed in Section~\ref{intuition} with the ones for $\mathfrak{a}^\mathfrak{A}(P)$, $\mathfrak{n}^\mathfrak{A}(P)$ and $\mathfrak{p}^\mathfrak{A}(P)$ displayed above it is easy to see that, when deleting the additional values, $\mathfrak{n}^\mathfrak{A}(P) = \mathfrak{a}^\mathfrak{A}(P) \setminus \mathfrak{c}^\mathfrak{A}(P)$ and $\mathfrak{p}^\mathfrak{A}(P) = \mathfrak{a}^\mathfrak{A}(P) \cup \mathfrak{c}^\mathfrak{A}(P)$ for every predicate symbol $P$.} From an axiomatic point of view, {\bf Dm} is obtained from {\bf Tm} by replacing \axT\ by an attenuated version of this axiom, namely, the well-known deontic axiom:\\

$
\begin{array}{ll}
\axD & \Box \alpha \rightarrow \Diamond \alpha
\end{array}$

\

From this, first-order extensions ${\bf Dm}^*$ and ${\bf Km}^*$ of {\bf Dm} and {\bf Km} can be defined, in the same way as ${\bf Tm}^*$ was obtained. Indeed, the definition of the corresponding Hilbert calculi is obvious: it is enough to add to the Hilbert calculi introduced in Part~I for these logics the axioms and the inference rule for the universal quantifier considered for ${\bf Tm}^*$. With respect to semantics,  let $V_k$ be the set of $k$ truth-values considered above, for $k=6,8$. Then, it is enough to expand the corresponding multialgebras with suitable multioperators $\tilde{Q}_k: (\mathcal{P}(V_k)-\{\emptyset\}) \to (\mathcal{P}(V_k)-\{\emptyset\})$ for every quantifier $Q\in \{\forall, \exists\}$ and $k=6,8$ extending the corresponding multifunctions over the four-element domain $V_k$ defined in Subsection~\ref{Nmatrix}. From the eight-valued multifuncions $\tilde{\wedge}$ and  $\tilde{\vee}$ for conjunction and disjunction displayed above it is straightforward to define, respectively,  $\tilde{\forall}_k$  and $\tilde{\exists}_k$ for $k=6,8$ in an analogous way as it was done in the four-valued case (we invite the reader to complete the details). Moreover, in order to validate \axNBF, the deterministic version $\tilde{Q}_k^d$ of these quantifiers can be easily defined by removing the occurrences of $F^-$ or $T^+$ from $\tilde{Q}_k(X)$ in the instances of $X$ where $\tilde{Q}_k(X)$ is not a singleton. Namely, for $k=6,8$: $\tilde{\forall}_k^d(X):=\{C^-\}$ for $X=\{C^+,C^-\}$ or $X=\{C^+,C^-,T^+\}$, while  $\tilde{\exists}_k^d(Y):=\{C^+\}$ for $Y=\{C^+,C^-\}$ or $Y=\{C^+,C^-,F^-\}$. We can extend this approach by considering the systems ${\bf D4m}^*$, ${\bf D45m}^*$, ${\bf K4m}^*$ and ${\bf K45m}^*$ obtained from the corresponding propositional systems studied in Part~I. Then, soundness and completeness of these systems with respect to the first-order Nmatrix semantics is easily obtained by adapting and extending the corresponding proofs for ${\bf Tm}^*$ detailed in the previous sections. Of course every $n$-ary predicate symbol $P$ should be interpreted in a six-valued or eight-valued first-order structure $\mathfrak{A}$ as a triple $P^\mathfrak{A} := (\mathfrak{a}^\mathfrak{A}(P),\mathfrak{n}^\mathfrak{A}(P),\mathfrak{p}^\mathfrak{A}(P))$ of subsets of $U^n$.\footnote{As a matter of fact, the observation made in Remark~\ref{equiv-sem} for ${\bf Tm}^*$ can be applied to these six-valued and eight-valued modal systems. Namely, interpreting a $n$-ary predicate symbol $P$ by means of a triple $P^\mathfrak{A}= (\mathfrak{a}^\mathfrak{A}(P),\mathfrak{n}^\mathfrak{A}(P),\mathfrak{p}^\mathfrak{A}(P))$ is equivalent to consider functions $P_\mathfrak{A}:U^n \to V_k$ for $k=6,8$. The latter approach is analogous to the non-deterministic first-order structures for paraconsistent logics based on Nmatrices considered in~\cite{CFG19}.} Once again, the details of these constructions are left to the interested reader.

None of the above mentioned six-valued and eight-valued systems validate axiom \axT. Invalidating \axT\ is not just a merely formal curiosity. We know that \axBF\ and \axBFe\ collapse in reflexive Kripkean models. In {\bf Tm}$^*$,  \axBFe\ implies \axBF\ and that seems to be a natural property, since strict implication is stronger than material implication. But that collapse is a direct consequence of axiom  \axT. The same consideration can be done with respect to \axCBF\ and \axCBFe. In any Ivlev-like first-order system in which \axT\ does not hold, the formulas \axBF\ and \axBFe\ will be totally independent, and the same holds for \axCBF\ and \axCBFe.

\

\section{Final Remarks} \label{final}

In this article, our investigations on non-normal modal logics  with finite-valued non-deterministic matrix semantics developed in Part~I (see~\cite{con:cer:per:19}) were extended to first-order languages.
By simplicity, and in order to fix the main ideas and the many advantages that this proposal can offer to the subject of first-order modal logics, just the quantified version of system {\bf Tm}, called ${\bf Tm}^*$, was analyzed in detail. The extension of this approach to other modal systems based on {\bf T4m}, {\bf T45m},  {\bf Dm}, {\bf D4m}, {\bf D45m}, {\bf Km}, {\bf K4m} or {\bf K45m} should  not be too difficult, as outlined in Section~\ref{other-sys}.

Concerning the formulas \axBF\ and \axCBF, we have proposed here a semantics in which these two formulas, at least at the first level, collapse (see Subsection~\ref{sectbarcan}). But this is not a consequence of the semantic clauses of atomic and propositional formulas. In fact, if our clauses regarding the universal operator were different, we could invalidate any of these formulas, as in the case of axiom \axNBF\ (recall Subsection~\ref{de_re}). As mentioned there, it should be interesting considering another notion of universal quantifier in which \axPBF\ would be blocked as well.

It worth noting that the way adopted in this paper defining the deterministic versions of the quantifiers intends to be the closest formal counterpart to our language intuitions involving terms like ``necessarily'' and ``for all''. 
This is not to say that we cannot discover another formal versions as much (or even more) intuitive than those presented here. We have not ruled out this possibility yet.

It is well-known that first-order classical logic is not decidable. But certain fragments of this logic are, in turn, decidable. An interesting case is its monadic fragment, that is, the fragment in which every  predicate is  unary and no function symbols are allowed. By contrast, Kripke has shown that almost all monadic first-order modal systems are undecidable.\footnote{See~\cite{kri:62}.} But this result does not automatically apply to {\bf Tm}$^*$ or any monadic fragment of the systems treated here. In fact, all propositional Ivlev-like systems proposed in Part~I are decidable. It would be interesting to know  how many of their monadic extensions are decidable.

Even though {\bf Tm}$^*$ is a system that does not compromise with the Kripkean thesis that proper names are rigid designators while definite descriptions are non-rigid designators, it would be interesting to consider a four-valued  non-deterministic modal semantics in which this occurs. A possible way would be extending the interpretation mappings  $\mathfrak{a}$ and $\mathfrak{c}$ to function symbols, while the individual constants would receive a unique interpretation in a given structure $\mathfrak{A}$. Besides, individual variables  would receive two kind of interpretations by means of the assignments, the actual one and the contingent one. We are aware, however, of the technical difficulties of this possible solution.

First-order logic has difficulties for dealing with  non-existent object that has a name, like ``Pegasus''. That happens because to each constant in the language we associate an individual in the domain. Since the existential quantifier traverses individuals in the domain of structure, in classical logic being named is equivalent to existing. Thus, if Pegasus is a winged horse, then Pegasus is a name of an individual of the domain and, thus, Pegasus exists. But we know that Pegasus is a winged horse that does not exist. In order to block this inference, we should use the so called Free Logics. It is a relevant philosophical problem to decide whether Pegasus is a being that does not exist but \emph{possibly} exists. Free Logics have made a very fruitful contribution to these investigations, especially when combined with semantics of possible worlds.\footnote{See \cite{sep-logic-free}.} We have no reason to be skeptical of the success of addressing these problems combining non-deterministic semantics with Free Logics.

\bibliographystyle{plain}
%\bibliography{BiblioConFarPer}

\end{document}